\documentclass{article}[11pt]
\usepackage{amsmath}
\usepackage{amsthm}
\usepackage{mathtools}
\usepackage{amsfonts}
\usepackage{subfigure}
\usepackage{caption}
\usepackage{tikz} 
\usepackage{amscd}

\usepackage{lipsum}
\usepackage{authblk}
\usepackage{fancyhdr}

\renewenvironment{abstract}{%
\hfill\begin{minipage}{0.95\textwidth}
\rule{\textwidth}{0.7pt}}
{ \rule{\textwidth}{1pt}\end{minipage}}

\makeatletter
\renewcommand\@maketitle{%
\hfill
\begin{minipage}{0.95\textwidth}
\vskip 2em
\let\footnote\thanks 
{\LARGE \@title \par }
\vskip 1.5em
{\large \@author \par}
\end{minipage}
\vskip 1em \par
}
\makeatother

\newtheorem{theorem}{Theorem}
\newtheorem{definition}{Definition}
\newtheorem{lemma}{Lemma}
\newtheorem{corollary}{Corollary}
\newtheorem{prop}{Proposition}

\begin{document}\Large

\title{BI-FILTRATION AND STABILITY OF TDA MAPPER FOR POINT CLOUD DATA}
\author[1]{WAKO BUNGULA}
\author[2]{ISABEL DARCY}
\affil[1]{UNIVERSITY OF WISCONSIN - LA CROSSE}
\affil[2]{UNIVERSITY OF IOWA}

\maketitle

\begin{abstract}
\normalsize
\textbf{Abstract} 

Carlsson, Singh and Memoli's TDA mapper takes a point cloud dataset and outputs a graph that depends on several parameter choices. Dey, Memoli, and Wang developed Multiscale Mapper for abstract topological spaces so that parameter choices can be analyzed via persistent homology. However, when applied to actual data, one does not always obtain filtrations of mapper graphs. DBSCAN, one of the most common clustering algorithms used in the TDA mapper software, has two parameters, \textbf{$\epsilon$} and \textbf{MinPts}. If \textbf{MinPts = 1} then DBSCAN is equivalent to single linkage clustering with cutting height \textbf{$\epsilon$}. We show that if DBSCAN clustering is used with \textbf{MinPts $>$ 2}, a filtration of mapper graphs may not exist except in the absence of free-border points; but such filtrations exist if DBSCAN clustering is used with \textbf{MinPts = 1} or \textbf{2} as the cover size increases, \textbf{$\epsilon$} increases, and/or \textbf{MinPts} decreases. However, the 1-dimensional filtration is unstable. If one adds noise to a data set so that each data point has been perturbed by a distance at most \textbf{$\delta$}, the persistent homology of the mapper graph of the perturbed data set can be significantly different from that of the original data set. We show that we can obtain stability by increasing both the cover size and \textbf{$\epsilon$} at the same time. In particular, we show that the bi-filtrations of the homology groups with respect to cover size and $\epsilon$ between these two datasets are \textbf{2$\delta$}-interleaved.

\end{abstract}

\section{Introduction}

TDA uses topological tools to analyze datasets of different types: topological spaces \cite{deybook}\cite{dey}, point clouds \cite{DBLP}, biological data \cite{Nicolau7265}, image data \cite{map_int}, graphs \cite{2018MOG}, and more. It provides insights to topological and geometric features of the underlying space of the data being analyzed. Persistence homology captures topological features that persists through time (or any other parameter). For a give parameter, one can obtain topological information, namely the $k^{th}$ homology group of the data, and the evolution of the topology as the parameter varies can be presented using a \textit{persistence diagram}.

  One of the main interests in TDA is to study stability of persistence diagrams. Suppose a given data, $X$, is perturbed by a small value $\delta$; denote the perturbed data $X_{\delta}$. The authors of \cite{Cohen-Steiner2007} and \cite{comptop} proved that the bottleneck distance between their two persistence diagrams is small and bounded. For 1-dimensional filtration (i.e. one varying parameter), it has been shown that the bottleneck distance can be used to measure how far (or close) two persistence diagrams are. However, for a multi-dimensional filtration (i.e. two or more varying parameters), there is no persistence diagram with nice presentation of the evolution of the topology \cite{DBLP}. Hence, the bottleneck distance can not be used. Chazal et. al., in \cite{chazel}, defined $\epsilon$-interleavings of persistence modules to introduce an idea of distance between persistence modules. This idea of distance between persistence modules makes it possible for stability claims. Michael Lesnick in \cite{DBLP} showed in the language of category theory that  $\xi$-interleaving is a nice generalization of the Bottleneck distance in that two $\epsilon$-interleaved persistence modules are algebraically similar.  
 
  If the data at hand is large and complex, then TDA may not be enough to analyze the topology of the underlying space. An algorithm, called \textit{TDA mapper}, was introduced by Gunnar Carlsson, Gurjeet Singh and Facundo Memoli, and it is used for simplification and visualization of complex and high-dimensional data while capturing topological features of the data \cite{map_int}. Studying the stability of mapper graphs has been of great interest among researchers  \cite{deybook}\cite{dey}\cite{stabMapp1}. 
 
  In \cite{deybook}\cite{dey}, Tamal Dey, Facundo Memoli and Yusu Wang proved stability results for \textit{Multiscale Mapper} when the input data is a topological space. Given a topological space $X$ equipped with a continuous function $f: X \rightarrow \mathbb{R}$, a filtration of covers of $\mathbb{R}$ gives rise to a filtration of covers of $X$ via the pre-image of $f$ which in turn gives rise to a filtration of the nerve of covers (or simplicial complexes, or abstract graphs, or mapper graphs). Furthermore, a filtration of simplicial complexes induces a filtration of homology groups \cite{deybook}\cite{dey}\cite{comptop}. For a topological space $X$ and its $\delta$-perturbation, $X_{\delta}$, Tamal Dey, Facundo Memoli and Yusu Wang (in  \cite{deybook}\cite{dey}) proved a stability theorem for the respective persistence diagrams given a cover of the image $f$ that satisfies a $(c,s)$-\textit{good} condition. 
 
  In \cite{stabMapp1}, Mathieu Carrière and Steve Oudot discuss, from  theoretical point of view, the relationship of the mapper graph and one of the Reeb graphs in order to predict the features of the mapper graph that is present (or not present) with respect to a given filter function, $f$, and the cover, $\mathbb{U}$, of the image of $f$. With this theoretical framework, the degree of stability of the features of the mapper graph can be quantified, and the convergence of the mapper graph to the Reeb graph can be guaranteed as the size of the elements of $\mathbb{U}$ goes to zero.
 
  While the focus of Dey et. al. is stability of \textit{Multiscale Mapper} when the input data is a topological space and the focus of Mathieu Carrière and Steve Oudot is a theoretical framework for stability of a one-dimensional mapper in relation to Reeb graph, the focus of this paper is bi-filtration and stability of mapper graphs when the input data is point cloud. Hence, what this paper presents is practical when dealing with the application of mapper algorithm and the input data in point cloud. 
 
  One of the results shown in this paper is that under certain conditions DBSCAN clustering and single-linkage clustering give a filtration of covers, hence stability of mapper graphs hold. The rest of this paper is outlined as follows. In section \ref{TDAPPLN}, the TDA mapper algorithm, whose input is a point cloud data and whose output is a simplicial complex, is presented. Different parameters of the TDA mapper yield different simplicial complexes, and studying multi-scale mapper helps us understand the relationship between different simplicial complexes. Section \ref{FSC_sec} revisits  \cite{deybook},\cite{dey} and \cite{comptop} that a filtration of covers induces a filtration of simplicial complexes that in turn induces a filtration of homology groups, and section  \ref{FMG} discuss how a filtration of mapper graphs can be realized with respect to varying TDA mapper parameters. Since clusters are defined based on clustering algorithms (\cite{CandC}), we investigate various clustering algorithms and if they give rise to a filtration of covers/clusters that gives rise to a filtration of mapper graphs. Section \ref{filtCL} gives a counter example to show that complete-linkage or average-linkage does not give a filtration of covers with respect to the parameter $bin$ size. Section \ref{dbsec} introduces DBSCAN clustering algorithm, and sections \ref{filtDB}, \ref{filtepsInc}, and \ref{filtmindec} discuss how DBSCAN gives a filtration of covers with respect to $bin$ size, $\epsilon$, and $MinPts$ respectively. 
 
  Section \ref{bifilt} discusses the bi-filtration of covers, simplicial complexes, and homology groups with respect to $bin$ size and $\epsilon$, and section \ref{interBiFit} reviews and reformulates interleaving of two bi-filtrations. Finally, section \ref{stabBiFilt} discuss the stability results, that is the bi-filtration of simplicial complexes and homology groups of two datasets $X$ and its $\delta$ perturbation $X_{\delta}$ are $2\delta$-interleaved.

\section{TDA Mapper Pipeline}\label{TDAPPLN}
Software implementing TDA mapper are written in various computer languages. AYASDI is a company in California based on TDA mapper that works with government, medical researchers, and financial institutes in analyzing data. Python-Mapper \cite{filfun}, Keppler-Mapper \cite{kmapper}, Mapper Interactive \cite{MapInt2011}, and Giotto-tda \cite{giotto-tda2020} are open source software written in Python.  Whereas TDA mapper \cite{tdaMappeR} is an open source software written in R.

\noindent The TDA mapper pipeline is as follows:
\begin{itemize}
\item The input is a point cloud data $X$ equipped with a filter function $f: X \rightarrow Z$, where often $Z = \mathbb{R}$. The choice of a filter function is made depending on the data and the researchers interest of study. 
\item The image of $f$ is equipped with a finite cover $\mathbb{U}$, where if $Z = \mathbb{R}$, $\mathbb{U}$ is a collection of intervals. The number of intervals and a percent overlap for the intervals are set by the user, which determines the length of the intervals.
\item The pre-image $f^{-1}(\mathbb{U})$ is a finite cover of $X$, called the pullback cover.  The elements of the pullback cover will be referred to as bins (or  overlapping bins due to nonempty intersections).
\item A clustering algorithm is applied in each of the overlapping bins.
\item Each cluster, $C_i$, is represented by a vertex, and a non-empty intersection between pairs of clusters is represented by an edge. Hence, the output is a graph. Note the clusters also form a cover of $X$ which will be called the cluster cover.  In the case where higher dimensional simplices (solid triangle, solid tetrahedron, ...) are formed via the nerve of the cluster cover, the output is a simplicial complex called the \v{C}ech complex.
\end{itemize}
\par
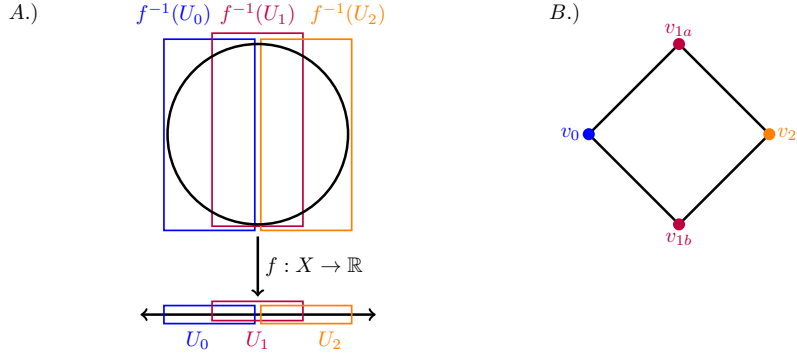
\begin{figure*}
\begin{center}
\scalebox{0.8}{
\begin{tikzpicture}

\draw[blue, thick] (-2.56,-1.6) rectangle (-1.05,1.58);

\draw[purple, thick] (-1.76,-1.53) rectangle (-0.25,1.68);

\draw[orange, thick] (-0.95,-1.6) rectangle (0.56,1.58);

\draw[color=black, very thick](-1,0) circle (1.5);

\draw[<->][black, very thick] (-2.95,-3) -- (0.96,-3);
\draw[->][black, very thick] (-1,-1.7) -- (-1,-2.7);
\fill (-1,-2.2)  circle (0) node[text=black, right]{$f:X\to\mathbb{R}$};

\draw[blue, thick] (-2.56,-3.15) rectangle (-1.05,-2.85);
\draw[purple, thick] (-1.76,-3.1) rectangle (-0.25,-2.78);
\draw[orange, thick] (-0.95,-3.15) rectangle (0.56,-2.85);

\fill (-2.,-3.15)  circle (0) node[text=blue, below]{$U_0$};
\fill (-1,-3.15)  circle (0) node[text=purple, below]{$U_1$};
\fill (0.2,-3.15)  circle (0) node[text=orange, below]{$U_2$};

\fill (-2.4,1.68)  circle (0) node[text=blue, above]{$f^{-1}(U_0)$};
\fill (-1,1.68)  circle (0) node[text=purple, above]{$f^{-1}(U_1)$};
\fill (0.5,1.68)  circle (0) node[text=orange, above]{$f^{-1}(U_2)$};

\draw[-][black, very thick] (4.5,0) -- (6,1.5);
\draw[-][black, very thick] (4.5,0) -- (6,-1.5);
\draw[-][black, very thick] (6,1.5) -- (7.5,0.);
\draw[-][black, very thick] (6,-1.5) -- (7.5,0.);

\fill [color=blue, very thick](4.5,0)  circle (0.1) node[text=blue, left]{$v_0$};
\fill [color=purple, very thick](6,1.5)  circle (0.1) node[text=purple, above]{$v_{1a}$};
\fill [color=purple, very thick](6,-1.5)  circle (0.1) node[text=purple, below]{$v_{1b}$};
\fill [color=orange, very thick](7.5,0.)  circle (0.1) node[text=orange, right]{$v_2$};

\fill [color=black, very thick](-4.5,2)  circle (0.) node[text=black, left]{$A.)$};
\fill [color=black, very thick](4.5,2)  circle (0.) node[text=black, left]{$B.)$};

\end{tikzpicture}}
\end{center}
\caption{A.) The input data, $X$, is a circle which is mapped to $\mathbb{R}$ via $f$. The image of $f$ is covered with $\mathbb{U} = \{U_0, U_1, U_2\}$, and the pre-image of $\mathbb{U}$, $f^{-1}(\mathbb{U}) = \{f^{-1}(U_0) , f^{-1}(U_1) , f^{-1}(U_2) \}$ covers $X$. B.) There is one connected component in each $f^{-1}(U_0)$ and $f^{-1}(U_2)$ represented by the two vertices $v_0$ and $v_2$. There are two connected components in $f^{-1}(U_1)$ represented by the two vertices $v_{1a}$ and $v_{1b}$.}
\label{TDAM_demo}	

\end{figure*}

 \noindent In figure \ref{TDAM_demo}.A, $X$ represents data in a circular shape equipped with the filter function $f:X\to\mathbb{R}$, $f(x, y) = x$. The image of $f$ is covered by a collection of three open intervals, $\mathbb{U}=\{U_0,U_1, U_2\}$. The pre-image of $\mathbb{U}$, $f^{-1}(\mathbb{U})=$$\{f^{-1}(U_0),$$  f^{-1}(U_1),  f^{-1}(U_2)\}$, covers $X$. If single-linkage clustering is used to identify connected components of the underlying circle in each bin,
there is one cluster each, $C_0$ and $C_1$, in $f^{-1}(U_0)$ and $f^{-1}(U_2)$, and there are two clusters, $C_{1a}$ and $C_{1b}$, in $f^{-1}(U_1)$. The cluster $C_i$ is represented by the vertex $v_i$. Figure \ref{TDAM_demo}.B shows the mapper graph where the vertices represent the clusters, and the edges represent nonempty intersection between pairs of clusters.

\section{Filtrations via Nerves of a Cover}\label{FSC_sec}
\noindent 

For a given topological space $X$ equipped with a continuous function $f: X \rightarrow \mathbb{R}$, Tamal Dey, Facundo Memoli, and Yusu Wang \cite{deybook}\cite{dey} showed that a filtration of covers of the image of $f$ induces a filtration of covers of $X$ by taking the pullback of $f$. They also showed that the filtration of covers of $X$ induces a filtration of the nerves of these covers, and by applying the homology functor, a filtration of homology groups of these nerves is obtained. 


\begin{figure}[h!]
\centering{
\begin{tikzpicture}

\node[draw,align=center] at (1, 1.)   (A) {Filtration of\\ Covers};
\node[draw,align=center]  at (5, 1)   (B) {Filtration of\\ Simplicial Complexes};
\node[draw,align=center]  at (10, 1)   (C) {Filtration of\\ Homology Groups};

\draw[->,line width=0.5mm](2.1,1)--(3.2,1);
\draw[->,line width=0.5mm](6.8,1)--(8.4,1);

\end{tikzpicture}}
\end{figure}

Below we state the formal definitions and theorems for completeness.

\par
\begin{definition}\label{nerveofU}(Nerve of a cover $\mathbb{U}$)
Given a cover $\mathbb{U} = \{U_\alpha\}_{\alpha \in A}$ of a topological space X, the nerve of the cover $\mathbb{U}$ is the simplicial complex $\mathbb{N}(\mathbb{U})$, called the \v{C}ech complex, whose vertex set is the index set A, and where a subset $\{\alpha_0, \alpha_1, \alpha_2, \dots, \alpha_n\} \subseteq A$ spans a $n$-simplex in $\mathbb{N}(\mathbb{U})$ if and only if $U_{\alpha_0} \cap U_{\alpha_1} \cap U_{\alpha_2} \cap \dots \cap U_{\alpha_n} \neq \emptyset$. 
\end{definition}
  The Nerve lemma asserts that given a good cover (i.e. sets in the cover are convex) of a topological space, $X$, the nerve of the cover and $X$ are homotopy equivalent. We can use a filtration of the nerves of the covers to study the topological evolution of $X$. Specifically, a filtration of of covers induces a filtration of simplicial complexes via their nerves, which induces a filtration of homological groups.

\begin{definition}  \cite{deybook}\cite{dey} \label{filtcover}(Filtration of covers)
Given a family of covers, $\mathbb{U} = \{\mathbb{U}_{\lambda}\}$ equipped with a family of maps $\{u^{\lambda_i, \lambda_j} : \mathbb{U}_{\lambda_i} \rightarrow \mathbb{U}_{\lambda_j}, \forall \lambda_i \leq \lambda_j\}$ for some parameter $\lambda$ where $u^{\lambda_i, \lambda_j} : \mathbb{U}_{\lambda_i} \rightarrow \mathbb{U}_{\lambda_j}$ such that if $u^{\lambda_i, \lambda_j}(U) = V$ then $U \subseteq V$ (i.e. $\mathbb{U}_{\lambda_i} \leq \mathbb{U}_{\lambda_j}$), we say there is a filtration of covers if 
$u^{\lambda_j, \lambda_k} \circ u^{\lambda_i, \lambda_j}= u^{\lambda_{i}, \lambda_k}$. 
\end{definition}
  Let $X$ be a topological space with $N$ coverings $\{\mathbb{U}^i\}_{i=1}^N = \{\{U^i_{\alpha}\}_{\alpha \in A_i}\}_{i=1}^N$. Also, let $h$ be a family of functions $\dots \rightarrow A_i \xrightarrow{h^{i,j}} A_j \xrightarrow{h^{j,k}} A_k \rightarrow \dots$ where $i \leq j \leq k$ such that for $\alpha \in A_i$, $h^{i,j}(\alpha) = \beta$ implies $U^i_{\alpha} \subseteq U^j_{\beta}$, and $h^{j,k} \circ h^{i,j} = h^{i,k}$. Then $h$ induces a well-defined family of simplicial maps $\dots \rightarrow \mathbb{N}(\mathbb{U}^i) \xrightarrow{\mathbb{N}(h^{i,j})} \mathbb{N}(\mathbb{U}^j) \xrightarrow{\mathbb{N}(h^{j,k})} \mathbb{N}(\mathbb{U}^k) \rightarrow \dots$ satisfying $\mathbb{N}(h^{j,k}) \circ \mathbb{N}(h^{i,j}) =\mathbb{N}(h^{i,k})$. Here $\mathbb{N}(h^{i,j})$ is a simplicial map defined on the vertex set of $\mathbb{N}(\mathbb{U}^i)$ such that if $\{v_0, v_1, ..., v_n\}$ is a simplex in $\mathbb{N}(\mathbb{U}^i)$ then $\{\mathbb{N}(h^{i,j})(v_0), \mathbb{N}(h^{i,j})(v_1), ..., \mathbb{N}(h^{i,j})(v_n)\}$ is a simplex in $\mathbb{N}(\mathbb{U}^j)$ where $\mathbb{N}(h^{i,j})(v_k) = v_{h^{i,j}(k)}$. Note an $m$-simplex may be mapped to an $n$-simplex such that $m\leq n$ since $h$ need not be one-to-one.
 
  In general, we have a family of covers of a topological space $X$ equipped with a parameter $\lambda$, and we define a filtration of covers as follows.

\begin{definition} \label{filtSimpComp}(Filtration of simplicial complexes)  
Given a family of simplicial complexes equipped with a family of simplicial maps $\{\phi^{\lambda_i, \lambda_j} : \mathbb{N}(\mathbb{U}_{\lambda_i}) \rightarrow \mathbb{N}(\mathbb{U}_{\lambda_j}), \forall \lambda_i \leq \lambda_j\}$ for some parameter $\lambda$, we say there is a filtration of simplicial complexes if 
$\phi^{\lambda_j, \lambda_k} \circ \phi^{\lambda_i, \lambda_j}= \phi^{\lambda_{i}, \lambda_k}$. That is, we have the following sequence of simplicial maps. 
\begin{center}
\begin{tikzpicture}

	\filldraw [black] (-1.715, 6.5) node[anchor=north] {$\mathbb{N}(\mathbb{U}_{\lambda_i})$};
	\filldraw [black] (0.985, 6.5) node[anchor=north] {$\mathbb{N}(\mathbb{U}_{\lambda_j})$};
	\filldraw [black] (3.715, 6.5) node[anchor=north] {$\mathbb{N}(\mathbb{U}_{\lambda_k})$};

	\draw[thick,->] (-1.0115, 6.2) -- (0.2, 6.2) node[anchor= south east]{$\phi^{\lambda_i, \lambda_j}$};
	\draw[thick,->] (1.685, 6.2) -- (3.0, 6.2) node[anchor= south east]{$\phi^{\lambda_j, \lambda_k}$};

	\draw[thick,->] (4.415, 6.2) -- (5., 6.2) ;
	\draw[thick,->] (-3.3, 6.2) -- (-2.6, 6.2) ;

	\node (a) at (-1.715, 6.5) {};
	\node (b) at (0.985, 6.5) {};
	\node (c) at (3.715, 6.5) {};
	\draw[->] (a)  to [in=90, out=90, looseness=0.5] (c);
	\filldraw [black] (1.485, 7.5) node[anchor=south east] {$\phi^{\lambda_i, \lambda_k}$};

\end{tikzpicture}
\end{center}

\end{definition}

\begin{theorem}\cite{munk}\label{filt_cover_to_filt_simpcompx}
Suppose there is a filtration of covers $\{u^{\lambda_i, \lambda_j} : \mathbb{U}_{\lambda_i} \rightarrow \mathbb{U}_{\lambda_j}, \forall \lambda_i \leq \lambda_j\}$, then for a fixed $n \in \mathbb{Z}$, this filtration of covers induces a well-defined family of simplicial maps $\{\phi^{\lambda_i, \lambda_j} : \mathbb{N}(\mathbb{U}_{\lambda_i}) \rightarrow \mathbb{N}(\mathbb{U}_{\lambda_j}), \forall \lambda_i \leq \lambda_j\}$ such that $\phi_n^{\lambda_j, \lambda_k} \circ \phi_n^{\lambda_i, \lambda_j} = \phi_n^{\lambda_i, \lambda_k}$.
\par
\noindent That is, a filtration of covers induces a filtration of simplicial complexes.
\end{theorem}
\noindent Theorem \ref{filt_cover_to_filt_simpcompx} states that given a filtration of cover, a filtration of simplicial complexes is obtained via the nerve of the cover. Theorem \ref{filt_homo} states that a filtration of simplicial complexes induces a filtration of homology groups. 

\begin{theorem}\cite{comptop}\label{filt_homo}
Suppose there is a filtration of simplicial complexes \\
$\{\phi^{\lambda_i, \lambda_j} : \mathbb{N}(\mathbb{U}_{\lambda_i}) \rightarrow \mathbb{N}(\mathbb{U}_{\lambda_j}), \forall \lambda_i \leq \lambda_j\}$, then for fixed $n \in \mathbb{Z}$, this filtration of simplicial complexes induce a well-defined family of homomorphisms,$\{f_n^{\lambda_i, \lambda_j}: H_n(\mathbb{N}(\mathbb{U}_{\lambda_i})) \rightarrow H_n(\mathbb{N}(\mathbb{U}_{\lambda_j}) ),  \forall \lambda_i \leq \lambda_j\}$, between the respective homology groups such that $f_n^{\lambda_j, \lambda_k} \circ f_n^{\lambda_i, \lambda_j} = f_n^{\lambda_i, \lambda_k}$.
\par
\noindent That is, a filtration of simplicial complexes induces a \textit{filtration} of homology groups.
\end{theorem}

\subsection{Filtration of Mapper Graphs}\label{FMG}
 
\par \noindent In TDA mapper, there are several parameters, such as interval length (\textit{bin} size) and percent overlap, chosen by the user. If two different sets of parameter choices are made, two different mapper graph outputs may be obtained. One or the other (or both) may contain topological information about the underlying space of the dataset. The question is ``How does one obtain the most accurate topological information about the underlying space of the dataset?" To answer this question, we study \textit{Multiscale Mapper} \cite{deybook}\cite{dey} applied to point cloud data; that is, we shall study the mapper graph outputs for a sequence of parameter choices. The idea is to analyze what topological feature(s) persist through these parameter choices.  

Recall that to create the mapper graph (or simplicial complex) from a data set $X$ using a filter function $f:X \rightarrow Z$, we first create a cover of $Z$ and pullback this cover to create a cover of $X$.  Each set in the pullback cover is then clustered.  These clusters form another cover of $X$ which we will call the cluster cover.  If we have a filtration of covers of $Z$, then we also have a filtration of the pullback covers of $X$.  However, as demonstrated in sections \ref{filtCL} and \ref{dbsec}, we do not always get a filtration of cluster covers from a filtration of pullback covers.  If we do have a filtration of cluster covers, then by \cite{deybook}\cite{dey} we have a filtration of their nerves and thus their homology groups per the previous section.  Thus for point cloud data $X$, we only need to check the first step in the sequence below:  Does the clustering algorithm result in a filtration of cluster covers given a filtration of covers of $X$

\begin{figure}[h!]
\centering{
\begin{tikzpicture}

\node[draw,align=center] at (1, 1.)   (A) {Filtration of\\ Covers of X};
\node[draw,align=center]  at (5.8, 1)   (B) {Filtration of\\ Cluster Covers};
\node[draw,align=center]  at (10.3, 1)   (B) {Filtration of\\ Mapper Graphs};
\node[draw,align=center]  at (15, 1)   (C) {Filtration of\\ Homology Groups};

\draw[->,line width=0.5mm](2.3,1)--(4.5,1)node[anchor= south east]{${\scalebox{1.5}{?}}\hskip 20pt$};
\draw[->,line width=0.5mm](7.2,1)--(8.9,1);
\draw[->,line width=0.5mm](11.8,1)--(13.4,1);

\end{tikzpicture}}
\end{figure}

  For a fixed set of parameters, the TDA mapper algorithm takes in a point cloud dataset, $X$, and outputs a mapper graph (or a simplicial complex). Recall from section \ref{TDAPPLN} the simplicial complex that the mapper output gives is obtained via the nerve of the clusters. An emphasis should be made that the mapper output is the simplicial complex obtained via the nerve of the clusters such that each vertex represents a cluster, each edge represents a non-empty intersection of two clusters, and an $n$-$simplex$ represents a non-empty intersection of $n + 1$ clusters. We say the nerve of clusters to mean the nerve of the cover because the set of clusters, $\{C_i\}$, is the cover of $X$. As the parameters of mapper vary, we prove the existence of a filtration of covers/clusters, which induces a filtration of simplicial complexes. In the following sections, we investigate whether or not various clustering algorithms give a filtration of covers/clusters via the nerve of clusters.


\subsection{No Filtration of Cluster Covers: Complete/Average-Linkage, $bin$ size.}\label{filtCL}

Single-linkage and DBSCAN are widely used because from topological point of view, these clustering algorithms give the correct connected components of a dataset. However, single-linkage will cluster two points that are far from each other as long as there is a chain of points connecting them, putting them in the same connected component.  Thus other clustering algorithms could be considered if closeness is a factor, but one would loose the topological advantages.  In addition to preserving topological properties like connectedness,
an additional advantage of using single-linkage and DBSCAN is that they give filtrations of Cech complexes/mapper graphs with respect to varying certain parameters per sections \ref{filtDB}, \ref{filtepsInc}, and \ref{filtmindec}, while most other clustering algorithms do not. 
TDA mapper in R only uses single-linkage as clustering algorithm, so we modified it in order to incorporate complete-linkage to illustrate the filtration issue with this and similar clustering algorithms.

\begin{figure*}[h]
\centering
	\centering
		A.) \includegraphics[width=6.5cm]{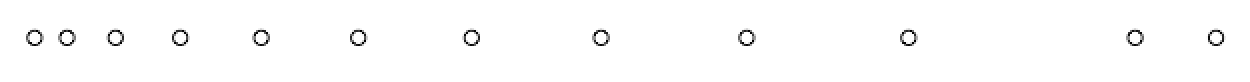}
	
		B.) \includegraphics[width=4.5cm]{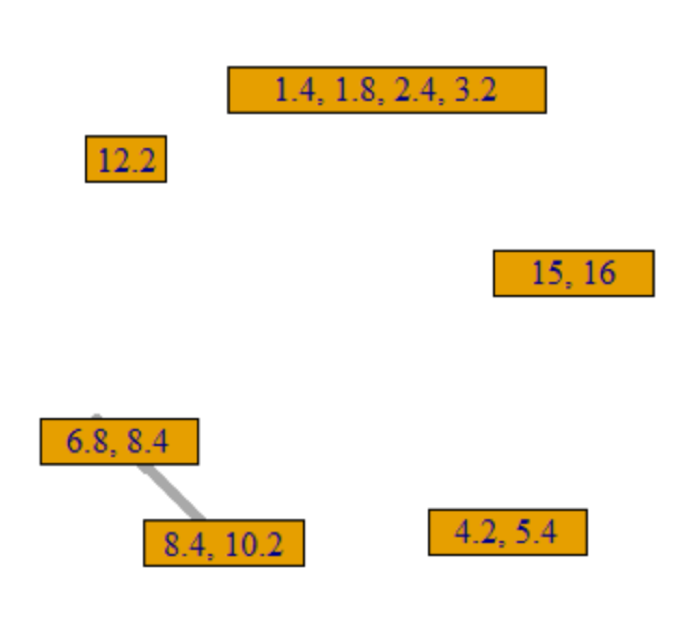}
		~~~~~
		C.) \includegraphics[width=4.5cm]{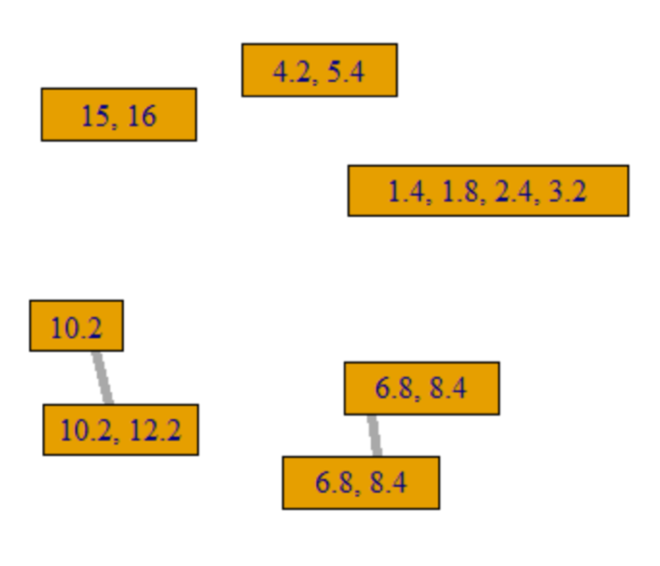}
	
\caption{A.) Dataset. B.) Mapper graph when complete-linkage is used with 2 intervals and $20\%$ overlap. C.) Mapper graph when complete-linkage is used with 2 intervals and $50\%$ overlap.}
\label{cl_dend1}

\end{figure*}

Suppose $X$ is the data, $X = \{1.4, 1.8, 2.4, 3.2, 4.2, 5.4, 6.8, 8.4, 10.2, 12.2, 15, 16\}$, containing 12 points (shown in figure \ref{cl_dend1}A). TDA mapper in R was run twice with identical parameters except for percent overlap. These parameter values are as follows with the TDA mapper R code given in parenthesis.
Filter function is $f:X\to \mathbb{R}$ such that $f(x)=x$ ($filter\_values = data[,1]$), number of intervals to cover the image of $f$ is 2 ($num\_intervals = 2$), and $num\_bins\_when\_clustering = 10$ which is used to determine cutting height when clustering. Let $I_1$ and $I_2$ be the two intervals with respect to $20\%$ overlap ($percent\_overlap = 20$), and $J_1$ and $J_2$ be the two intervals with respect to $50\%$ overlap ($percent\_overlap = 50$). The mapper graphs are shown in figures  \ref{cl_dend1}.B and \ref{cl_dend1}.C.
 
\par \noindent Note $f^{-1}(I_1)\subseteq f^{-1}(J_1)$ and $f^{-1}(I_2) \subseteq f^{-1}(J_2)$. Thus the covers form a filtration $\{f^{-1}(I_1), f^{-1}(I_2)\} \leq \{f^{-1}(J_1), f^{-1}(J_2)\}$.  We note that cluster $\{8.4, 10.2\}$ is in $f^{-1}(I_2)$ but is not contained in any of the clusters in $f^{-1}(J_2)$. This counter example shows that as $bin$ size increases (i.e. $f^{-1}(I_1)\subseteq f^{-1}(J_1)$ and $f^{-1}(I_2) \subseteq f^{-1}(J_2)$.), complete-linkage does not give filtration of cluster covers where the elements of the covers are the clusters obtained via complete-linkage. 
\noindent The same argument (with the same dataset and parameters) shows that average-linkage does not give filtration of cluster covers via clusters with respect to increasing $bin$ size.

\par
\par

\section{DBSCAN}\label{dbsec}

After discussing the DBSCAN algorithm, we will discuss the conditions under which DBSCAN gives filtration of cluster covers. Since single-linkage is a special case of DBSCAN, the results proven for DBSCAN also apply to single-linkage.

 The DBSCAN clustering method uses `density' to determine which points belong to which cluster \cite{DB}. The general idea is to cluster a set of dense points together whose underlying space is connected, and if there is set of low-density points, they are labeled as noise. Two parameters are needed to carry out the DBSCAN algorithm: these parameters are minimum number of points  $MinPts \in \mathbb{Z}$ and a radius $\epsilon \in \mathbb{R}$. The parameter $MinPts$ is the minimum number of points required in an $\epsilon$ neighborhood of a point $p$ for the $\epsilon$ neighborhood of $p$ to be considered a dense set. Note that Kepler-Mapper uses DBSCAN clustering, and the default values of the parameters are $\epsilon=0.5$ and $MinPts = 3$. An $\epsilon-neighborhood$ of a point $p$ is defined to be the set of points in the dataset that are at most $\epsilon$ distance away from $p$, and it is denoted by $N_{\epsilon}(p) = \{q \in X : dist(p,q) \leq \epsilon \}$. Note that $N_{\epsilon}(p)$ is a closed ball, and $X$ is the entire dataset.

\begin{figure*}[h]
\centering

\begin{tikzpicture}[scale=1.25]

\filldraw[black] (-1.0, -0.35) circle (0.9pt);
\filldraw[black] (-0.3, -0.235) circle (0.9pt);
\filldraw[blue] (-2., 0.0) circle (0.9pt)node[anchor=west] {$r$};
\filldraw[blue] (-1.8, 0.2) circle (0.9pt);
\filldraw[blue] (-2.3, -0.2) circle (0.9pt);
\draw [blue] (-2., 0.0)  circle (0.4cm);
\draw [blue] (-1.8, 0.2)  circle (0.4cm);
\draw [blue] (-2.3, -0.2)  circle (0.4cm);
\filldraw[black] (-0.8, 0.4) circle (0.9pt);
\filldraw[red] (-0.5, -0.15) circle (0.9pt); 
\filldraw[red] (-0.4, -0.05) circle (0.pt) node[anchor=east] {$p$};
\draw [red](-0.5, -0.19) circle (0.4cm);

\filldraw[black] (-0.65, 0.157) circle (0.9pt) ;
\filldraw[black] (-1.0, -0.025) circle(0.9pt);
\filldraw[black] (-0.6, -0.5) circle (0.9pt)node[anchor=west] {$q$};
\draw [black](-0.6, -0.5) circle (0.4cm);
\filldraw[black] (-0.4, 0.15) circle (0.9pt);

\filldraw[black] (-0.6, 0.3) circle (0.9pt);
\filldraw[black] (-0.44, 0.6) circle (0.9pt);
\filldraw[black] (-0.2, 0.2) circle (0.9pt);
\filldraw[black] (-0.25, 0.4) circle (0.9pt);
\filldraw[black] (-0.1, 0.5) circle (0.9pt);
\filldraw[black] (-0.05, 0.3) circle (0.9pt);
\filldraw[black] (-0.05, 0.8) circle (0.9pt);
\filldraw[black] (0.2, 0.5) circle (0.9pt);
\filldraw[black] (0.24, 0.7) circle (0.9pt);

\end{tikzpicture}

\caption{Let $MinPts = 5$ and $\epsilon$ be the radius of the circles centered at $p$, $q$, and $r$. Point $p$ is a core point because $|N_{\epsilon}(p)|$=$5 \geq MinPts$. Point $q$ is a border point because $q \in N_{\epsilon}(p)$ and $|N_{\epsilon}(q)|$=$3 < MinPts$. Point $r$ is neither a core point nor a border point, hence it is noise.} 
\label{ddr}
\end{figure*}
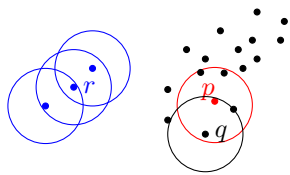

  Figure \ref{ddr} illustrates core points, border points, and noise points which are defined as follows. 

\begin{itemize}
\item A point $p$ is a \textit{core point} if $|N_{\epsilon}(p)| \geq MinPts$.
\item A point $q$ is a \textit{border point} if $|N_{\epsilon}(q)| < MinPts$, and $q \in N_{\epsilon}(p)$ for $p$ a core point.
\item A point $r$ is \textit{noise} if $r$ is neither a core point nor a border point.  \label{noisedef}
\end{itemize}
 
  Neither border points nor noise points exist when $MinPts = 1$. Every point in $X$ is a core point. In the case when $MinPts = 1$, clusters obtained by applying DBSCAN with respect to $\epsilon$ and $Minpts$ are the same cluster obtained by applying single-linkage with a cutting height $\epsilon$. Hence, single-linkage is a special case of DBSCAN.

  The following definitions, taken from \cite{DB}, are required to formally define a cluster obtained by the DBSCAN algorithm.

\begin{definition}{(directly density-reachable)} \cite{DB} \label{def_ddr}
A point $q$ is directly density-reachable from a point $p$ wrt. $\epsilon$ and $MinPts$ if $p$ is a core point and $q \in N_{\epsilon}(p)$.
\end{definition}

\begin{definition}{(density-reachable)}  \cite{DB} \label{def_dr}
A point $q$ is density-reachable from a point $p$ wrt. $\epsilon$ and $MinPts$ if there is a sequence of points $p=\alpha_1, \alpha_2, \alpha_3, ..., \alpha_{n-1}, \alpha_n = q$ such that $\alpha_{i+1}$ is directly density-reachable form $\alpha_i$. 
\end{definition}

\begin{definition}{(density-connected)}  \cite{DB} \label{def_dc}
A point $p$ is density-connected to a point $q$ wrt. $\epsilon$ and $MinPts$ if there is a point $o$ such that both, $p$ and $q$ are density-reachable from $o$ wrt. $\epsilon$ and $MinPts$.
\end{definition}

\begin{figure*}
\centering
	A.) \includegraphics[scale=0.3]{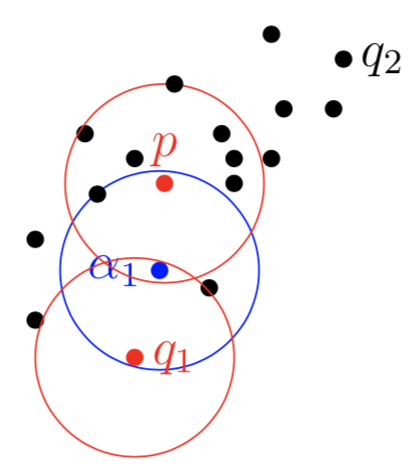}	
	~~~~~~~~~
	B.) \includegraphics[scale=0.3]{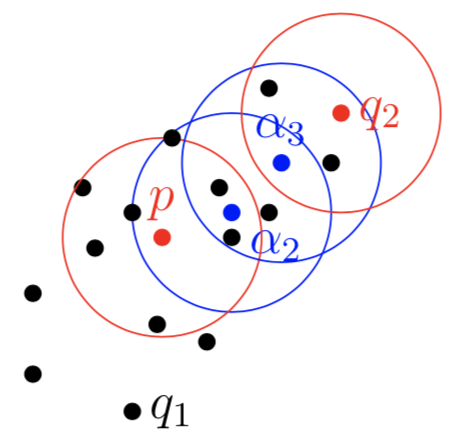}
	~~~~~~~~~
	C.) \includegraphics[scale=0.3]{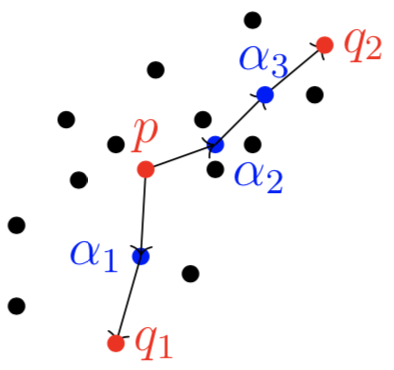}
\caption{Let $Minpts = 5$ and $\epsilon=$ the radius of the circles. A.) Point $q_1$ is directly density-reachable from $\alpha_1$ and density-reachable from $p$ wrt. $\epsilon$ and $MinPts$. B.) Point $q_2$ is directly density-reachable from $\alpha_3$ and density-reachable from $p$ wrt. $\epsilon$ and $MinPts$. C.) Points $q_1$ and $q_2$ are density-connected wrt. $\epsilon$ and $MinPts$.}

\label{DB_DDRC}	
\end{figure*}
\noindent Figures \ref{DB_DDRC} A and B show that points $q_1$ and $q_2$ are density-reachable from point $p$ wrt. $\epsilon$ and $MinPts$, and thus points $q_1$ and $q_2$ are density-connected wrt. $\epsilon$ and $MinPts$ (per figure \ref{DB_DDRC}.C).

\begin{definition}(cluster)\label{correctDBSCAN}
Let a dataset $X = C_1 \sqcup C_2 \sqcup ...\sqcup C_n \sqcup N$. That is, we are partitioning $X$ into disjoint sets such that $N$ is a set of noise points and $C_i$ is a cluster wrt. $\epsilon$ and $MinPts$ satisfying:
\begin{itemize}
\item (Maximality) $\forall p, q$: if $p \in C_i$, $q \notin \sqcup_{j=1}^{i-1}{C_j}$, and $q$ is density-reachable from $p$ wrt. $\epsilon$ and MinPts, then $q \in C_i$.
\item (Connectivity) $\forall p, q \in C_i: p$ is density-connected to $q$ wrt $\epsilon$ and MinPts. 
\end{itemize}
\end{definition}

\par \noindent \textit{Remark-1} Definition \ref{correctDBSCAN} is a modification of definition 5 in \cite{DB} taking into consideration that order of the dataset matters as illustrated below.
 
\par \noindent \textit{Remark-2} The memory complexity for the current DBSCAN clustering algorithm written in python is $\mathcal{O}(n\times d)$ where $d$ is the average number of points in an $\epsilon$ neighborhood of core points and $n$ is $MinPts$ \cite{complexity}.

\begin{figure*}[h]
\centering
	A.) \includegraphics[scale=0.2]{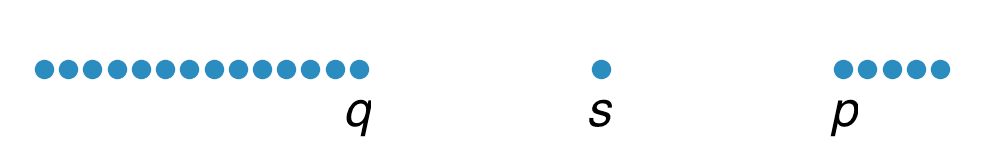}

\vspace{0.5cm}

	B.) \includegraphics[scale=0.2]{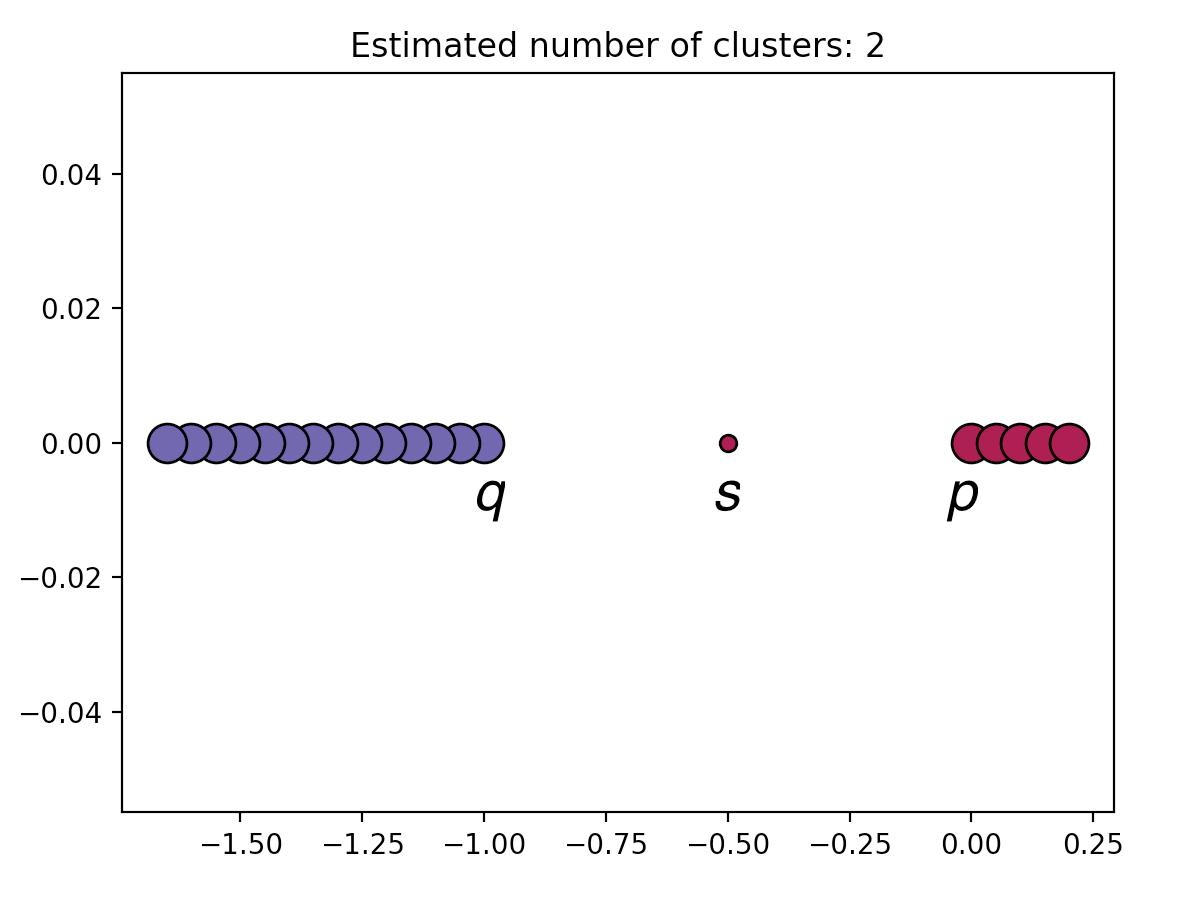}
	~~~~~~~~~~~~~~
	C.) \includegraphics[scale=0.2]{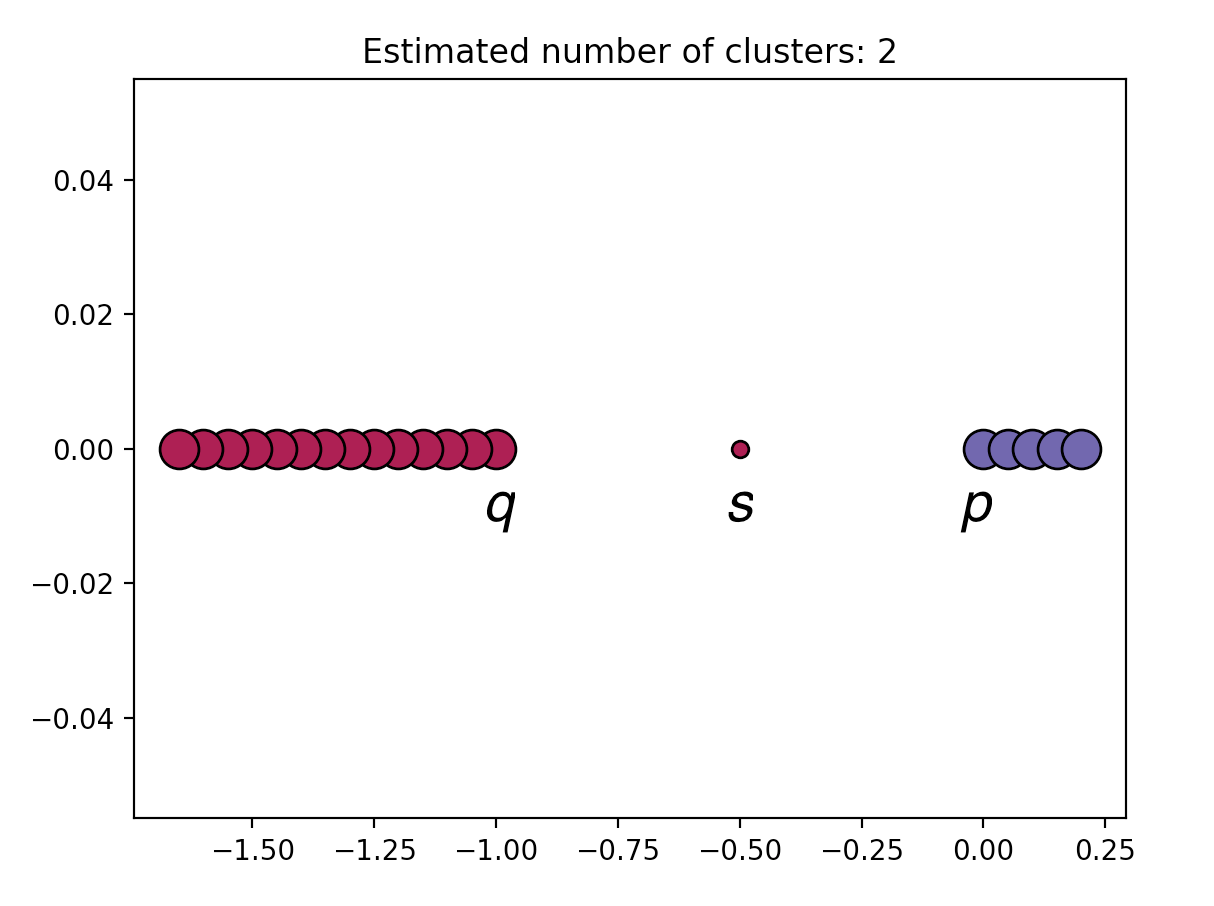}
	\label{dbstability2}
\caption{Point $s$ is a \textit{free-border} point since it belongs to different clusters depending on the order of points in the dataset. A.) Dataset. B.) Point $s$ is clustered with point $p$ if $p$ or any of the point to the right of $p$ is listed first. C.) Point $s$ is clustered with point $q$ if $q$ or any of the point to the left of $q$ is listed first. }
\label{dbstability}
\end{figure*}

\noindent
We will now illustrate how the output of DBSCAN depends on how the data set is ordered.
Consider the dataset shown in figure \ref{dbstability}.A. Let $d(p,s) = d(q,s) = 0.5$, and the distance between any two consecutive points to the left of point $q$ and any two consecutive points to the right of point $p$ be less than or equal to $0.16$. Point $s$ is directly density-reachable from both point $q$ and point $p$ wrt $\epsilon = 0.5$ and $Minpts = 5$, but neither point $p$ nor point $q$ are directly density reachable from point $s$. Thus, 
\begin{itemize}
\item $s \in N_{\epsilon}(p)$, $s \in N_{\epsilon}(q)$, $|N_{\epsilon}(p)| \geq 5$, and  $|N_{\epsilon}(q)| \geq 5$. 
\item $N_{\epsilon}(s) < 5$; i.e. $s$ is a border point. And, $p$ is not density-connected to $q$.
\end{itemize}

\noindent There are two clusters wrt $\epsilon$ and $MinPts$, one cluster containing point $p$ and the points to the right of $p$ and another cluster containing point $q$ and the points to the left of $q$. Since clusters do not intersect, the question is: ``To which cluster does point $s$ belong?" The answer lies in the order in which the data is listed. Let $X = \{x_1, x_2, ...\}$ represent the dataset in figure \ref{dbstability}.A. 
Point $s$ is clustered with $p$ if $x_1$ is $p$ or any of the points to the right of $p$ or if $x_1$ is $s$ and $x_2$ is $p$ or any of the points to the right of $p$. See figure \ref{dbstability}.B. Similarly, $s$ is clustered with $q$ if $x_1$ is $q$ or any of the points to the left of $q$ or if $x_1$ is $s$ and $x_2$ is $q$ or any of the points to the left of $q$. See figure \ref{dbstability}.C.

  The DBSCAN clustering algorithm takes as an input an ordered set of points, and two parameters, $\epsilon$ and $MinPts$. After the core points are identified wrt. $\epsilon$ and $MinPts$, clusters are formed wrt. $\epsilon$ and $MinPts$ based on their order. Any point that does not belong to any of the clusters is labeled as noise. 

\begin{definition}{(free-border point)}\label{fbp}
Let $s$ be a point in a given dataset. We say $s$ is a free-border point wrt. $\epsilon$ and $MinPts$ if there exist two points $p$ and  $q$ satisfying
\begin{itemize}
\item $s \in N_{\epsilon}(p)$ and  $|N_{\epsilon}(p)| \geq MinPts$,
\item $s \in N_{\epsilon}(q)$ and  $|N_{\epsilon}(q)| \geq MinPts$,
\item $|N_{\epsilon}(s)| < MinPts$, and 
\item $p$ is not density connected to $q$.
\end{itemize}
\end{definition}

\begin{lemma}{\label{fbpIFF}}
A point $s$ is a free-border point with respect to $\epsilon$ and $MinPts$ if and only if there exist two orderings of the dataset such that $s$ belongs to different clusters depending on the ordering.
\end{lemma}

\noindent Lemmas \ref{correct1} and \ref{correct2} are modified versions of Lemmas 1 and 2 of \cite{DB} respectively taking order into consideration.  

\begin{lemma}\label{correct1}
Let $X$ be a dataset.
\begin{itemize}
\item Let $p$ be a point in $X$ and $|N_{\epsilon}(p)| \geq MinPts.$ Then there exists a cluster C wrt. $\epsilon$ and $MinPts$ containing the point $p$. Moreover, if $q$ is in C, then $q$ is density-reachable from $p$.

\item Let C be a cluster wrt. $\epsilon$ and $MinPts$, then there exist a point $p$ in C with $|N_{\epsilon}(p)| \geq MinPts.$ That is $p$ is core point of $C$ wrt. $\epsilon$ and $MinPts$. If $q$ is in C, then $q$ is density-reachable from $p$.
\end{itemize}
\end{lemma}

\noindent In short, lemma \ref{correct1} declares that a cluster is determined by any of its core points, and the core points of a cluster will remain core points of the same cluster no matter the order of the dataset. 

\begin{lemma}{\label{correct2}}
Let C be a cluster wrt. $\epsilon$ and $MinPts$, and let $C_{core}= \{p$ $|$ $p$ is a core point of $C \}$. That is if $p \in C_{core}$, $|N_{\epsilon}(p)| \geq MinPts$. Then $\forall p \in C_{core}$, $C_{core} = \{q$ $|$ $q$ is a core point, and $q$ is a density-reachable from  $p \}$.
  
\end{lemma}
 
\noindent Lemma \ref{correct2} states the following. Suppose two core points, $p_1$ and $p_2$, are in the same cluster ( i.e. they represent the same cluster). Then $p_1$ and $p_2$ are in the same cluster regardless of the order of the dataset.

In the sections below, we will show that if there are no free border points, a filtration of covers of $X$ will give a filtration of cluster covers of $X$ as bin size increases (section \ref{filtDB}), the $\epsilon$ parameter of DBSCAN increases (section \ref{filtepsInc}), or the MinPts parameter of DBSCAN decreases (section \ref{filtmindec}).

\subsection{Filtration of cluster covers: DBSCAN, $\mathbb{B} = bin$ size.}\label{filtDB}

Recall that DBSCAN takes two parameters $\epsilon$ and $MinPts$. We will investigate three cases when DBSCAN is applied to a dataset: $MinPts = 1$,  $MinPts = 2$, and  $MinPts > 2$.  If $MinPts = 1$, the clusters obtained wrt. $\epsilon$ and $MinPts = 1$ are the same clusters if single-linkage is applied with a cutting height $h = \epsilon$. If DBSCAN is used to cluster a dataset $X$ with the parameters $\epsilon$ and $MinPts = 1$, then every point is a core point. Because every data point in $X$ is a core point, there are no free-border points. In fact, there are neither border points nor noise points wrt. $\epsilon$ and $MinPts=1$. Hence a filtration of cluster covers, as $bin$ size increases, is obtained when DBSCAN is used wrt. $\epsilon$ and $MinPts = 1$, and it is the same filtration of cluster covers when single-linkage is used where the cutting height is $\epsilon$ for all $bin$. 
 
  If $MinPts = 2$, by corollary \ref{DBMP=2} (see below), there are no free-border points, and thus there is a filtration of cluster covers as $bin$ size increases. If $MinPts > 2$,  filtration of cluster covers is not guaranteed because there can exist free-border points. Although, if we assume no free-border points exist wrt. $\epsilon$ and $MinPts$, there is filtration of cluster covers as $bin$ size increases as stated in lemma \ref{filcovDB}.
 
  The R version of DBSCAN has a parameter \textit{borderPoints} whose default argument is \textit{TRUE}. Though only free-border points are the problem, one can choose to set all border points as noise points by setting parameter \textit{borderPoints = FALSE}, in which case data points are partitioned into core points and noise points. That is, a cluster only contains core points. Moreover, a filtration of cluster covers is obtained if one chooses to not include border points in a cluster, but one should note that the filtration of simplicial complexes obtained by ignoring border points may be different from the filtration of simplicial complexes obtained by ignoring just free-border points since an intersection between two clusters could consist of only border points. However the filtration at the level of vertices would be the same.
 
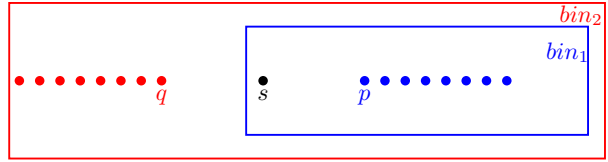
\begin{figure*}
\begin{center}
\scalebox{0.9}{
\begin{tikzpicture}

\draw[red, thick] (-2.75,-3.15) rectangle (6.05,-0.85);
\draw[blue, thick] (0.75,-2.8) rectangle (5.8,-1.2);

\fill [color=red, thick](-2.6,-2.)  circle (0.07);
\fill [color=red, thick](-2.3,-2.)  circle (0.07);
\fill [color=red, thick](-2.0,-2.)  circle (0.07);
\fill [color=red, thick](-1.7,-2.)  circle (0.07);
\fill [color=red, thick](-1.4,-2.)  circle (0.07);
\fill [color=red, thick](-1.1,-2.)  circle (0.07);
\fill [color=red, thick](-0.8,-2.)  circle (0.07);
\fill [color=red, thick](-0.5,-2.)  circle (0.07) node[text=red, below]{$q$};

\fill [color=black, thick](1.,-2.)  circle (0.07) node[text=black, below]{$s$};

\fill [color=blue, thick](2.5,-2.)  circle (0.07) node[text=blue, below]{$p$};
\fill [color=blue, thick](2.8,-2.)  circle (0.07);
\fill [color=blue, thick](3.1,-2.)  circle (0.07);
\fill [color=blue, thick](3.4,-2.)  circle (0.07);
\fill [color=blue, thick](3.7,-2.)  circle (0.07);
\fill [color=blue, thick](4,-2.)  circle (0.07);
\fill [color=blue, thick](4.3,-2.)  circle (0.07);
\fill [color=blue, thick](4.6,-2.)  circle (0.07);

\fill [color=blue, thick](5.5,-1.3)  circle (0.0) node[text=blue, below]{$bin_1$};
\fill [color=blue, thick](5.7,-0.75)  circle (0.0) node[text=red, below]{$bin_2$};

\end{tikzpicture}}
\end{center}
\caption{Point $s$ is directly density reachable from both points $q$ and $p$ wrt. $\epsilon = 0.5$ and $MinPts = 5$. DBSCAN clusters $s$ and $p$ together in $bin1$. But, in $bin2$, $s$ is clustered with $q$ if any one of the red points is listed before the blue points.}
\label{stabilityofDBSCAN}
\end{figure*}

  How do free-border points fail to give a filtration of cluster covers? Consider the dataset, $X=\{q,s,p,...\}$, shown in figure \ref{stabilityofDBSCAN}. Note $q$ is the first point in the ordered dataset. Let $MinPts = 5$ and $\epsilon = dist(p,s) = dist(q,s)$. Note $bin1 \subseteq  bin2$. When DBSCAN is applied in $bin1$, one cluster, $C_p^{bin1}$, is formed wrt. $\epsilon$ and $MinPts$ containing $s$, $p$, and all the blue points to the right of $p$. When DBSCAN is applied to $bin2$, two clusters are formed wrt. $\epsilon$ and $MinPts$. One cluster, $C_q^{bin2}$, contains points $s$, $q$, and all the red points to the left of $q$ since $q$ is listed first in the ordered dataset. The second cluster, $C_p^{bin2}$, contains points $p$ and all the blue points to right of $p$. Although $bin1 \subseteq bin2$, $C_p^{bin1} \not\subseteq C_p^{bin2}$ since $s \in C_{p}^{bin1}$ but $s \notin C_{p}^{bin2}$. Hence, because free-border points could belong to different clusters in different bins, filtration of cluster covers does not always exist. 

  If there exist no free-border points in X wrt. $\epsilon$ and $MinPts$, then there is a filtration of cluster covers, where a cover of a dataset $X$ is the set of clusters. The following lemma states containment of clusters, implying filtration of cluster covers, as $bin$ size increases in the absence of free-border points. Recall lemma \ref{correct1} and lemma \ref{correct2} imply that a cluster is determined by any of core points, and we denote by $C_p^{bin1}$ a cluster determined by a core point $p$ in $bin1$.

\begin{lemma}\label{filcovDB}
Suppose there are no free-border points when DBSCAN is used to cluster a dataset, $X$. Let $ bin1 \subseteq bin2$. Denote by $C_{p}^{bin1}$  a cluster in $bin1$ determined by a core point $p$ wrt. $\epsilon$ and $MinPts$, and $C_{p}^{bin2}$ a cluster in $bin2$ determined by a core point $p$ wrt. $\epsilon$ and $MinPts$. Then, $C_{p}^{bin1} \subseteq C_{p}^{bin2}$. 
\end{lemma}
\par
\begin{proof}
Suppose $C_{p}^{bin1}$ is a cluster in $bin1$ and $q \in C_{p}^{bin1}$. 
Since $q$ is density-reachable from $p$ wrt. $\epsilon$ and $MinPts$ in $bin1$, $q$ is density-reachable from $p$ wrt. $\epsilon$ and $MinPts$ in $bin2$ because $bin1 \subseteq bin2$. Since there are no free-border points, no core point from which $q$ is density reachable is in $bin2$ that is not in $C_{p}^{bin2}$. Then $q \in C_{p}^{bin2}$. Thus, $C_{p}^{bin1} \subseteq C_{p}^{bin2}$. 

\end{proof}

  Lemma \ref{filcovDB} implies theorem \ref{filtrationofcoversasbinincreases} that states in the absence of free-border points, filtration of cluster covers is obtained. And, if filtration of cluster covers is realized, we get filtration of the nerve of covers/clusters and a filtration of homology groups  \cite{deybook}\cite{dey}. Let us denote the collection of overlapping bins of a given dataset by $\mathbb{B}_i = \{bin_i\}$, and we write $\mathbb{B}_i \leq \mathbb{B}_j$ if for all $bin_i \in \mathbb{B}_i$ there is $bin_j \in \mathbb{B}_j$ such that $bin_i \subseteq bin_j$.
\par

\begin{theorem}(Filtration of cluster covers in the parameter $\mathbb{B}$.)\label{filtrationofcoversasbinincreases}
Let $X$ be a dataset, and DBSCAN is used to cluster X. Let $\epsilon$ and $MinPts$ be fixed. If no free-border points exist wrt. $\epsilon$ and $MinPts$, there is a filtration of cluster covers, $\{c^{\mathbb{B}_i, \mathbb{B}_j} :\mathbb{C}_{\mathbb{B}_i} \rightarrow \mathbb{C}_{\mathbb{B}_j}, \forall \mathbb{B}_i \leq \mathbb{B}_j\}$ where $\mathbb{C}_{\mathbb{B}_i}$ is a cover of $X$ with respect to $\mathbb{B}_i$. 

\end{theorem}

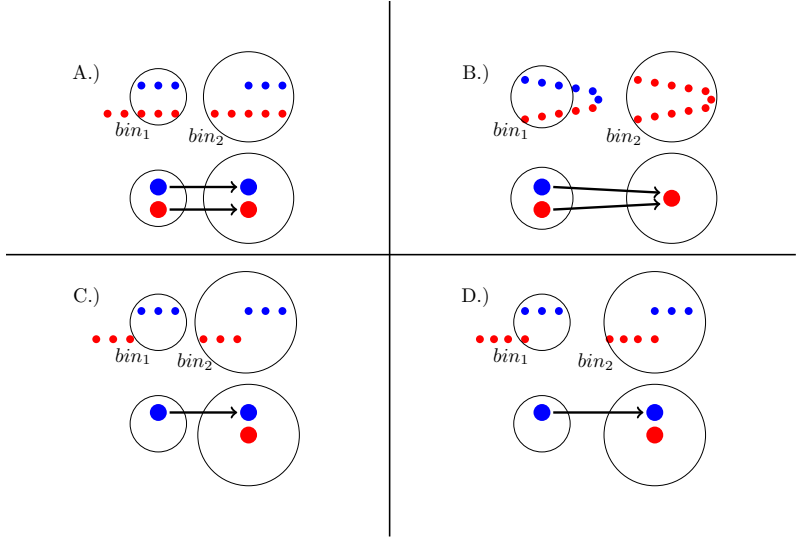
\begin{figure*}
\begin{center}
\scalebox{0.75}{

\begin{tikzpicture}

\fill (-4.2,0.2)  circle (0) node[text=black, left]{A.)};
\fill (2.7,0.2)  circle (0) node[text=black, left]{B.)};
\fill (-4.2,-3.75)  circle (0) node[text=black, left]{C.)};
\fill (2.7,-3.75)  circle (0) node[text=black, left]{D.)};

\draw[-][black, thick] (-6,-3.) -- (8,-3.);
\draw[-][black, thick] (0.8,1.5) -- (0.8,-8);

\fill (-3.3,-4.8)  circle (0) node[text=black, left]{$bin_1$};
\fill (-2.2,-4.9)  circle (0) node[text=black, left]{$bin_2$};

\fill [color=blue, thick](-3.6,0)  circle (0.07); 
\fill [color=blue, thick](-3.3,0)  circle (0.07); 
\fill [color=blue, thick](-3.,0)  circle (0.07); 

\fill [color=red, thick](-4.2,-0.5)  circle (0.07); 
\fill [color=red, thick](-3.9,-0.5)  circle (0.07); 
\fill [color=red, thick](-3.6,-0.5)  circle (0.07); 
\fill [color=red, thick](-3.3,-0.5)  circle (0.07); 
\fill [color=red, thick](-3.0,-0.5)  circle (0.07); 

\draw(-3.3,-0.2)  circle (0.5); 

\fill [color=blue, thick](-1.7,0)  circle (0.07); 
\fill [color=blue, thick](-1.4,0)  circle (0.07); 
\fill [color=blue, thick](-1.1,0)  circle (0.07); 

\fill [color=red, thick](-2.3,-0.5)  circle (0.07); 
\fill [color=red, thick](-2,-0.5)  circle (0.07); 
\fill [color=red, thick](-1.7,-0.5)  circle (0.07); 
\fill [color=red, thick](-1.4,-0.5)  circle (0.07); 
\fill [color=red, thick](-1.1,-0.5)  circle (0.07); 

\draw(-1.7,-0.2)  circle (0.8); 

\fill [color=blue, very thick](-3.3,-1.8)  circle (0.15); 
\fill [color=red, very thick](-3.3,-2.2)  circle (0.15); 
\fill [color=blue, very thick](-1.7,-1.8)  circle (0.15); 
\fill [color=red, very thick](-1.7,-2.2)  circle (0.15); 

\draw[->][black, very thick] (-3.1,-1.8) -- (-1.9,-1.8);
\draw[->][black, very thick] (-3.1,-2.2) -- (-1.9,-2.2);

\draw(-3.3,-2)  circle (0.5); 
\draw(-1.7,-2)  circle (0.8); 

\fill (3.4,-0.8)  circle (0) node[text=black, left]{$bin_1$};
\fill (5.4,-0.9)  circle (0) node[text=black, left]{$bin_2$};

\fill [color=blue, thick](3.2,0.1)  circle (0.07); 
\fill [color=blue, thick](3.5,0.05)  circle (0.07); 
\fill [color=blue, thick](3.8,0)  circle (0.07); 
\fill [color=blue, thick](4.1,-0.05)  circle (0.07); 
\fill [color=blue, thick](4.4,-0.1)  circle (0.07); 
\fill [color=blue, thick](4.5,-0.25)  circle (0.07); 

\fill [color=red, thick](3.2,-0.6)  circle (0.07); 
\fill [color=red, thick](3.5,-0.55)  circle (0.07); 
\fill [color=red, thick](3.8,-0.5)  circle (0.07); 
\fill [color=red, thick](4.1,-0.45)  circle (0.07); 
\fill [color=red, thick](4.4,-0.4)  circle (0.07); 

\draw(3.5,-0.2)  circle (0.55); 

\fill [color=red, thick](5.2,0.1)  circle (0.07); 
\fill [color=red, thick](5.5,0.05)  circle (0.07); 
\fill [color=red, thick](5.8,0)  circle (0.07); 
\fill [color=red, thick](6.1,-0.05)  circle (0.07); 
\fill [color=red, thick](6.4,-0.1)  circle (0.07); 
\fill [color=red, thick](6.5,-0.25)  circle (0.07); 

\fill [color=red, thick](5.2,-0.6)  circle (0.07); 
\fill [color=red, thick](5.5,-0.55)  circle (0.07); 
\fill [color=red, thick](5.8,-0.5)  circle (0.07); 
\fill [color=red, thick](6.1,-0.45)  circle (0.07); 
\fill [color=red, thick](6.4,-0.4)  circle (0.07); 

\draw(5.8,-0.2)  circle (0.8); 

\fill [color=blue, very thick](3.5,-1.8)  circle (0.15); 
\fill [color=red, very thick](3.5,-2.2)  circle (0.15); 
\fill [color=red, very thick](5.8,-2)  circle (0.15); 

\draw[->][black, very thick] (3.7,-1.8) -- (5.6,-1.9);
\draw[->][black, very thick] (3.7,-2.2) -- (5.6,-2.1);

\draw(3.5,-2)  circle (0.55); 
\draw(5.8,-2)  circle (0.8); 

\fill (-3.3,-0.8)  circle (0) node[text=black, left]{$bin_1$};
\fill (-2,-0.9)  circle (0) node[text=black, left]{$bin_2$};

\fill [color=blue, thick](-3.6,-4)  circle (0.07); 
\fill [color=blue, thick](-3.3,-4)  circle (0.07); 
\fill [color=blue, thick](-3.,-4)  circle (0.07); 

\fill [color=red, thick](-3.8,-4.5)  circle (0.07); 
\fill [color=red, thick](-4.1,-4.5)  circle (0.07); 
\fill [color=red, thick](-4.4,-4.5)  circle (0.07); 

\fill [color=blue, thick](-1.7,-4)  circle (0.07); 
\fill [color=blue, thick](-1.4,-4)  circle (0.07); 
\fill [color=blue, thick](-1.1,-4)  circle (0.07); 

\fill [color=red, thick](-1.9,-4.5)  circle (0.07); 
\fill [color=red, thick](-2.2,-4.5)  circle (0.07); 
\fill [color=red, thick](-2.5,-4.5)  circle (0.07); 

\draw(-3.3,-4.2)  circle (0.5); 
\draw(-1.75,-4.2)  circle (0.9); 

\fill [color=blue, very thick](-3.3,-5.8)  circle (0.15); 
\fill [color=blue, very thick](-1.7,-5.8)  circle (0.15); 
\fill [color=red, very thick](-1.7,-6.2)  circle (0.15); 

\draw[->][black, very thick] (-3.1,-5.8) -- (-1.9,-5.8);

\draw(-3.3,-6)  circle (0.5); 
\draw(-1.7,-6.2)  circle (0.9); 

\fill (3.4,-4.8)  circle (0) node[text=black, left]{$bin_1$};
\fill (4.9,-4.9)  circle (0) node[text=black, left]{$bin_2$};

\fill [color=blue, thick](3.2,-4)  circle (0.07); 
\fill [color=blue, thick](3.5,-4)  circle (0.07); 
\fill [color=blue, thick](3.8,-4)  circle (0.07); 

\fill [color=red, thick](3.2,-4.5)  circle (0.07); 
\fill [color=red, thick](2.9,-4.5)  circle (0.07); 
\fill [color=red, thick](2.65,-4.5)  circle (0.07); 
\fill [color=red, thick](2.4,-4.5)  circle (0.07); 

\fill [color=blue, thick](5.5,-4)  circle (0.07); 
\fill [color=blue, thick](5.8,-4)  circle (0.07); 
\fill [color=blue, thick](6.1,-4)  circle (0.07); 

\fill [color=red, thick](5.5,-4.5)  circle (0.07); 
\fill [color=red, thick](5.2,-4.5)  circle (0.07); 
\fill [color=red, thick](4.95,-4.5)  circle (0.07); 
\fill [color=red, thick](4.7,-4.5)  circle (0.07); 

\draw(3.5,-4.2)  circle (0.5); 
\draw(5.5,-4.2)  circle (0.9); 

\fill [color=blue, very thick](3.5,-5.8)  circle (0.15); 
\fill [color=blue, very thick](5.5,-5.8)  circle (0.15); 
\fill [color=red, very thick](5.5,-6.2)  circle (0.15); 

\draw[->][black, very thick] (3.7,-5.8) -- (5.3,-5.8);

\draw(3.5,-6)  circle (0.5); 
\draw(5.5,-6.2)  circle (0.9); 

\end{tikzpicture}}
\end{center}

\caption {The top row in each subfigure shows the dataset as $bin$ size increases while the bottom row shows the resulting filtration of mapper graphs. A.) As $bin$ size increases (i.e. $bin1 \subseteq bin2$), the blue cluster in $bin1$ is equal to the blue cluster in $bin2$, and the red cluster in $bin1$ is the subset of the red cluster in $bin2$. The blue vertex (representing the blue cluster) is mapped to the blue vertex, and the red vertex (representing the red vertex) is mapped to the red vertex. B.) As $bin$ size increases (i.e. $bin1 \subseteq bin2$), two clusters (red and blue) in $bin1$ merge to one red cluster in $bin2$. Hence, the two vertices representing the two clusters merge to one vertex. C.) As $bin$ size increases (i.e. $bin1 \subseteq bin2$), a cluster (blue) is born. The vertex representing the new cluster is also born. D.) As $bin$ size increases (i.e. $bin1 \subseteq bin2$), a cluster (blue) is born. The noise points in $bin1$ are included in the newborn cluster in $bin2$. The vertex representing the new cluster is also born in $bin2$.}
\label{filtcover4case}

\end{figure*}

\begin{corollary}\label{DBMP=2} Let DBSCAN be used to cluster X with respect to $\epsilon$ and $MinPts = 1, 2$. Then, there is filtration of cluster covers as $bin$ size increases.
\end{corollary}
\begin{proof}
By lemma \ref{filcovDB}, it only suffices to prove that there are no free-border points wrt. $\epsilon$ and $MinPts = 1, 2$.  Suppose $MinPts =1$ which implies every point is a core point, hence no free-border point exists. Suppose $MinPts=2$, and for the sake of contradiction, suppose that there is a free-border point wrt. $\epsilon$ and $MinPts = 2$. By definition \ref{fbp}, if $s$ is a free-border point, $|N_{\epsilon}(s)| < MinPts = 2$, and there is a core point $p$ such that $s \in N_{\epsilon}(p)$. This is a contradiction because $\{s, p\} \in N_{\epsilon}(s)$. Therefore, there are no free-border points wrt. $\epsilon$ and $MinPts = 2$.

\end{proof}
  An exhaustive case of what could happen as $bin$ size increases is shown in figures \ref{filtcover4case}. Note $bin$ size is a parameter of mapper, whereas $\epsilon$ and $MinPts$ are parameters of DBSCAN. In the absence of free-border points, we have proven the existence of filtration of cluster covers as the $bin$ size increases, hence by the results from  \cite{deybook}\cite{dey} (reviewed in section \ref{FSC_sec}), we get theorem \ref{filtbinTheo} as a consequence. That is, we have the existence of filtration of the nerve of cluster covers as the parameter $bin$ size increases.

\subsection{Filtration of cluster covers: DBSCAN, $\epsilon$.} \label{filtepsInc}
Suppose $\epsilon_0 \leq \epsilon_1$. If $p$ is a core point wrt. $\epsilon_0$ and $MinPts$, $|N_{\epsilon_0}(p)| \geq MinPts$. Since $\epsilon_0 \leq \epsilon_1$, $|N_{\epsilon_1}(p)| \geq MinPts$. Denote by $C_p^{\epsilon}$ a cluster determined by a core point $p$ wrt.  $\epsilon$ and $MinPts$ in $bin1 \in \mathbb{B}$. Therefore, $p$ is a core point determining a cluster $C_p^{\epsilon_0}$ wrt. $\epsilon_0$ and $MinPts$ in $bin1 \in \mathbb{B}$, and $p$ is also a core point determining a cluster $C_p^{\epsilon_1}$ wrt. $\epsilon_1$ and $MinPts$ in $bin1\in \mathbb{B}$.
\par
\begin{lemma}\label{filtepslem}
Suppose $X$ is a dataset, $\epsilon_0 \leq \epsilon_1$, and $MinPts$ and $\mathbb{B}$ is fixed. Let $C_p^{\epsilon_0}$ be a cluster wrt. $\epsilon_0$ and $MinPts$ and $C_p^{\epsilon_1}$ be a cluster wrt. $\epsilon_1$ and $MinPts$. If there are no free-border points wrt.  $\epsilon_1$ and $MinPts$, then $C_p^{\epsilon_0} \subseteq C_p^{\epsilon_1}$.

\end{lemma}
\begin{proof}
Suppose $q \in C_p^{\epsilon_0}$, then $q$ is density reachable from $p$ wrt. $\epsilon_0$ and $MinPts$.  Since $ \epsilon_0 \leq \epsilon_1$, then $q$ is density reachable from $p$ wrt. $\epsilon_1$ and $MinPts$. Because there are no free-border points wrt. $\epsilon_1$ and $MinPts$ and $p$ remains to be a core point of $C_p^{\epsilon_1}$ wrt. $\epsilon_1$ and $MinPts$, $q \in C_p^{\epsilon_1}$.

\end{proof}

\noindent Lemma \ref{filtepslem} implies theorem \ref{filtrationofcoversasepsincreases} that is there is a filtration of cluster covers as the DBSCAN parameter $\epsilon$ increases in the absence of free-border points.

\begin{theorem}(Filtration of cluster covers in the parameter $\epsilon$.)\label{filtrationofcoversasepsincreases}
Let $X$ be a dataset, and DBSCAN is used to cluster $X$. Let $\mathbb{B}$ and $MinPts$ be fixed. If no free-border points exist wrt. $\epsilon_i$ and $MinPts$ for all $i$, there is a filtration of cluster covers, $\{c^{\epsilon_i, \epsilon_j} : \mathbb{C}_{\epsilon_i} \rightarrow \mathbb{C}_{\epsilon_j}, \forall \epsilon_i \leq \epsilon_j\}$ where $\mathbb{C}_{\epsilon_i}$ is a cover of $X$ with respect to $\epsilon_i$. 

\end{theorem}

\subsection{Filtration of cluster covers: DBSCAN, Decreasing $MinPts$.} \label{filtmindec}
  Suppose $MinPts_0 \geq MinPts_1$. If $p$ is a core point wrt. $\epsilon$ and $MinPts_0$, then $|N_{\epsilon}(p)| \geq MinPts_0 \geq MinPts_1$. Therefore, if $p$ is a core point determining the cluster $C_p^{MinPts_0}$ wrt. $\epsilon$ and $MinPts_0$, then $p$ is also a core point determining the cluster $C_p^{MinPts_1}$ wrt. $\epsilon$ and $MinPts_1$.

\par
\begin{lemma}\label{filtminlem}
Suppose $X$ is a dataset, $MinPts_0 \geq MinPts_1$, and let $\epsilon$ and $\mathbb{B}$ be fixed. Let $C_p^{MinPts_0}$ be a cluster wrt. $\epsilon$ and $MinPts_0$ and $C_p^{MinPts_1}$ be a cluster wrt. $\epsilon$ and $MinPts_1$. If there are no free-border point wrt. $\epsilon$ and $MinPts_1$, then $C_p^{MinPts_0} \subseteq C_p^{MinPts_1}$. 

\end{lemma}

\begin{proof}
Suppose $q \in C_p^{MinPts_0}$, then $q$ is density-reachable from $p$ wrt. $\epsilon$ and $MinPts_0$. Since $MinPts_0 \geq MinPts_1$,  $q$ is density-reachable from $p$ wrt. $\epsilon$ and $MinPts_1$. Because there are no free-border points wrt. $\epsilon$ and $MinPts_1$ and $p$ is a core point of $C_p^{MinPts_1}$ wrt. $\epsilon$ and $MinPts_1$, $q \in C_p^{MinPts_0}$.

\end{proof}
\par

\begin{theorem}(Filtration of cluster covers in the parameter $MinPts$.)\label{filtrationofcoversasMinPtsincreases}
Let $X$ be a dataset, and DBSCAN is used to cluster $X$. Let $\mathbb{B}$ and $\epsilon$ be fixed. If no free-border points exist wrt. $MinPts_j$ and $\epsilon$, there is a filtration of clusters covers, $\{c^{MinPts_i, MinPts_j} : \mathbb{C}_{MinPts_i} \rightarrow \mathbb{C}_{MinPts_j}, \forall MinPts_i \geq MinPts_j\}$ where $\mathbb{C}_{MinPts_i}$ is a cover of $X$ with respect to $MinPts_i$. 

\end{theorem}

\noindent Lemma \ref{filtminlem} implies theorem \ref{filtrationofcoversasMinPtsincreases} that there is a filtration of cluster covers as the DBSCAN parameter $MinPts$ decreases in the absence of free-border points. 
\noindent In this section, we have shown that DBSCAN gives filtrations of covers of a dataset in three parameters, $\mathbb{B}$, $\epsilon$, and $MinPts$ in the absence of free-border points. In section \ref{filt_SC_and_HG}, we will state the existence of filtrations of simplicial complexes and homology groups.

\section{Filtration of Simplicial Complexes and Homology Groups}\label{filt_SC_and_HG}
Applying the result from  \cite{deybook}\cite{dey} as reviewed in section \ref{FSC_sec}, theorem \ref{filtrationofcoversasbinincreases} gives theorem \ref{filtbinTheo}, theorem \ref{filtrationofcoversasepsincreases} gives theorem \ref{filtEpsTheo}, and theorem \ref{filtrationofcoversasMinPtsincreases} gives theorem \ref{filtMinPtsTheo}. It follows, there are filtrations of the nerve of clusters induced by the filtrations of covers of X in the three parameters, $\mathbb{B}$, $\epsilon$, and $MinPts$. Moreover, there exist filtrations of homology groups induced by the filtrations of simplicial complexes.

\begin{theorem}(Filtration of simplicial complexes and homology groups in the parameter $\mathbb{B}$.){\label{filtbinTheo}}
Let $X$ be a dataset, and DBSCAN is used to cluster $X$. Let $k \in \mathbb{Z}$, $\epsilon$, and $MinPts$ be fixed. If no free-border points exist wrt. $\epsilon$ and $MinPts$, 
\begin{itemize}
\item There is a filtration of simplicial complexes, $\{\phi^{\mathbb{B}_i, \mathbb{B}_j} : \mathbb{N}(\mathbb{C}_{\mathbb{B}_i}) \rightarrow\mathbb{N}(\mathbb{C}_{\mathbb{B}_j}), \forall \mathbb{B}_i \leq \mathbb{B}_j\}$ where $\mathbb{N}(\mathbb{C}_{\mathbb{B}_i})$ is a simplicial complex via the nerve of covers (or clusters) whose vertex set represent clusters wrt. $\epsilon$ and $MinPts$. 
\item There is a filtration of homology groups, $\{f^{\mathbb{B}_i, \mathbb{B}_j} : H_k(\mathbb{N}(\mathbb{C}_{\mathbb{B}_i})) \rightarrow H_k(\mathbb{N}(\mathbb{C}_{\mathbb{B}_j})), \forall \mathbb{B}_i \leq \mathbb{B}_j\}$ where $H_k(\mathbb{N}(\mathbb{C}_{\mathbb{B}_i}))$ is the $k^{th}$ homology group of $X$ with respect to $\mathbb{B}_i$. 
\end{itemize}
\end{theorem}

\begin{theorem}(Filtration of simplicial complexes and homology groups in the parameter $\epsilon$){\label{filtEpsTheo}}
Let $X$ be a dataset, and DBSCAN is used to cluster $X$. Let $k$, $\mathbb{B}$, and $MinPts$ be fixed. For a sequence of $\epsilon$ values $\epsilon_0 \leq \epsilon_1 \leq \epsilon_2 \leq \dots \leq \epsilon_m$, let there be no free-border points wrt. $\epsilon_i$ and $MinPts$. 
\begin{itemize}
\item There is a filtration of simplicial complexes, $\{\phi^{\epsilon_i, \epsilon_j} : \mathbb{N}(\mathbb{C}_{\epsilon_i}) \rightarrow \mathbb{N}(\mathbb{C}_{\epsilon_j}), \forall \epsilon_i \leq \epsilon_j\}$ where $\mathbb{N}(\mathbb{C}_{\epsilon_i})$ is a simplicial complex via the nerve of clusters (or the mapper graph) whose vertex set represent clusters wrt. $\epsilon_i$ and $MinPts$.
\item There is a filtration of homology groups, $\{f^{\epsilon_i, \epsilon_j} : H_k(\mathbb{N}(\mathbb{C}_{\epsilon_i})) \rightarrow H_k(\mathbb{N}(\mathbb{C}_{\epsilon_j})), \forall \epsilon_i \leq \epsilon_j\}$ where $H_k(\mathbb{N}(\mathbb{C}_{\epsilon_i}))$ is the $k^{th}$ homology group of $X$ with respect to $\epsilon_i$. 

\end{itemize}

\end{theorem}

\begin{theorem}(Filtration of simplicial complexes and homology groups in the parameter $MinPts$){\label{filtMinPtsTheo}}
Let $X$ be a dataset, and DBSCAN is used to cluster $X$. Let $k$, $\mathbb{B}$, and $\epsilon$ be fixed. For a sequence of $MinPts$ values $MinPts_0 \geq MinPts_1 \geq MinPts_2 \geq \dots \geq MinPts_m$, let there be no free-border points wrt. $\epsilon$ and $MinPts_i$. 
\begin{itemize}
\item There is a filtration of simplicial complexes, $\{\phi^{MinPts_i, MinPts_j} : \mathbb{N}(\mathbb{C})_{MinPts_i} \rightarrow \mathbb{N}(\mathbb{C})_{MinPts_j}, \forall MinPts_i \geq MinPts_j\}$ where $\mathbb{N}(\mathbb{C})_{MinPts_i}$ is a simplicial complex via the nerve of clusters (or the mapper graph) whose vertex set represent clusters wrt. $\epsilon$ and $MinPts_i$.
\item There is a filtration of homology groups, $\{f^{MinPts_i, MinPts_j} : H_k(\mathbb{N}(\mathbb{C}))_{MinPts_i} \rightarrow H_k(\mathbb{N}(\mathbb{C}))_{MinPts_j}, \forall MinPts_i \geq MinPts_j\}$ where $H_k(\mathbb{N}(\mathbb{C}_{\epsilon_i}))$ is the $k^{th}$ homology group of $X$ with respect to $\epsilon_i$. 

\end{itemize}

\end{theorem}

\section{Bi-Filtrations and Stability}\label{bifilt}

\begin{figure}
\begin{center}

\scalebox{0.65}{
\begin{tikzpicture}

\fill (-4.8,0.2)  circle (0) node[text=black, left]{A.)};
\fill (0.7,0.2)  circle (0) node[text=black, left]{B.)};
\fill (5.2,0.2)  circle (0) node[text=black, left]{C.)};

\fill [color=blue, very thick](-4.5,0)  circle (0.1);
\fill [color=blue, very thick](-4.2,0)  circle (0.1); 
\fill [color=blue, very thick](-3.9,0)  circle (0.1); 
\fill [color=blue, very thick](-3.6,0)  circle (0.1); 
\fill [color=blue, very thick](-3.3,0)  circle (0.1); 
\fill [color=blue, very thick](-3.0,0)  circle (0.1); 
\fill [color=blue, very thick](-2.7,0)  circle (0.1); 
\fill [color=blue, very thick](-2.4,0)  circle (0.1); 
\fill [color=blue, very thick](-2.1,0)  circle (0.1); 

\fill [color=red, very thick](-4.2,-1)  circle (0.1); 
\fill [color=red, very thick](-3.9,-1)  circle (0.1); 
\fill [color=red, very thick](-3.6,-1)  circle (0.1); 
\fill [color=red, very thick](-3.3,-1)  circle (0.1); 
\fill [color=red, very thick](-3.0,-1)  circle (0.1); 
\fill [color=red, very thick](-2.7,-1)  circle (0.1); 
\fill [color=red, very thick](-2.4,-1)  circle (0.1); 

\draw(-3.3,-0.2)  circle (1.3); 
\draw(2.2,-0.1)  circle (1.29); 
\draw(6.9,-0.2)  circle (1.7); 
\fill (-4.3,-1.2) circle (0) node[text=black, left]{$bin$};
\fill (1.2,-1.1) circle (0) node[text=black, left]{$bin$};
\fill (5.7,-1.4) circle (0) node[text=black, left]{$bin_{\delta}$};

\fill [color=blue, very thick](1,0.5)  circle (0.1);
\fill [color=blue, very thick](1.3,0.5)  circle (0.1); 
\fill [color=blue, very thick](1.6,0.5)  circle (0.1); 

\fill [color=black, very thick](1.9,0)  circle (0.1); 
\fill [color=black, very thick](2.2,0)  circle (0.1); 
\fill [color=black, very thick](2.5,0)  circle (0.1); 

\fill [color=orange, very thick](2.8,0.5)  circle (0.1); 
\fill [color=orange, very thick](3.1,0.5)  circle (0.1); 
\fill [color=orange, very thick](3.4,0.5)  circle (0.1); 

\fill [color=red, very thick](1.3,-1.5)  circle (0.1); 
\fill [color=red, very thick](1.6,-1.5)  circle (0.1); 
\fill [color=red, very thick](1.9,-1.5)  circle (0.1); 
\fill [color=red, very thick](2.2,-1.5)  circle (0.1); 
\fill [color=red, very thick](2.5,-1.5)  circle (0.1); 
\fill [color=red, very thick](2.8,-1.5)  circle (0.1); 
\fill [color=red, very thick](3.1,-1.5)  circle (0.1); 

\fill [color=blue, very thick](5.7,0.5)  circle (0.1);
\fill [color=blue, very thick](6,0.5)  circle (0.1); 
\fill [color=blue, very thick](6.3,0.5)  circle (0.1); 

\fill [color=black, very thick](6.6,0)  circle (0.1); 
\fill [color=black, very thick](6.9,0)  circle (0.1); 
\fill [color=black, very thick](7.2,0)  circle (0.1); 

\fill [color=orange, very thick](7.5,0.5)  circle (0.1); 
\fill [color=orange, very thick](7.8,0.5)  circle (0.1); 
\fill [color=orange, very thick](8.1,0.5)  circle (0.1); 

\fill [color=red, very thick](6,-1.5)  circle (0.1); 
\fill [color=red, very thick](6.3,-1.5)  circle (0.1); 
\fill [color=red, very thick](6.6,-1.5)  circle (0.1); 
\fill [color=red, very thick](6.9,-1.5)  circle (0.1); 
\fill [color=red, very thick](7.2,-1.5)  circle (0.1); 
\fill [color=red, very thick](7.5,-1.5)  circle (0.1); 
\fill [color=red, very thick](7.8,-1.5)  circle (0.1); 

\end{tikzpicture}}
\end{center}

\caption{DBSCAN is applied to cluster wrt. $\epsilon$ = distance between two adjacent points and $MinPts$ = 3. A.) There are two clusters wrt. $\epsilon$ and $MinPts$ in $bin$. B.) There are points in $X_{\delta}$ not in $bin$. C.) There are four clusters wrt. $\epsilon$ and $MinPts$ in $bin_\delta$, but there are two cluster wrt. $\epsilon + 2\delta$ and $MinPts$.}
\label{stblty12}

\end{figure}
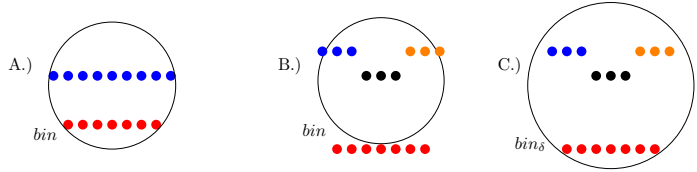

We would motivate why stability needs to be discussed by demonstrating DBSCAN is not stable under small perturbation. Given a dataset $X$, suppose we perturb $X$ by at most $\delta$, and denote the perturbed dataset by $X_{\delta}$. That is, there is a function $\Delta: X\rightarrow X_{\delta}$ such that $dist(x, \Delta(x)) \leq \delta$ where $x\in X$. The distance between $X$ and $X_{\delta}$ is $dist(X, X_{\delta}) = max\{min_{x \in X}\{dist(x,y) : x \in X, y \in X_{\delta}\}\} \leq \delta$. Applying DBSCAN to cluster $X$ and $X_{\delta}$ could yield significant difference in topology of the dataset. So, a \textit{small change in the dataset yields a significant change in the topology}. Consider the example shown in figure \ref{stblty12} where $X$ is shown in figure \ref{stblty12}.A and $X_{\delta}$ is shown in figure \ref{stblty12}.C. If we set $\epsilon$ to the distance between two adjacent points of $X$ and $MinPts = 3$, then there are two clusters wrt. $\epsilon$ and $MinPts$. The first cluster, $C_1^{bin}$, contains nine points, and the second cluster, $C_2^{bin}$, contains seven points. Applying DBSCAN to $X_{\delta}$ gives four clusters, where three clusters, $C_1^{bin_{\delta}}$, $C_2^{bin_{\delta}}$, $C_3^{bin_{\delta}}$, contain three points each and one cluster, $C_4^{bin_{\delta}}$, contains seven points. We will state the conditions under which this issue of instability is resolved in section \ref{stabBiFilt}. 

\begin{figure}
\centering
\resizebox{1.95in}{1.5in}{

\begin{tikzpicture}
  \node (o) at (0,0) {$\mathbb{C}_{(\mathbb{B}_0, \epsilon_0)}$};
  
  \node (a) at (0,2) {$\mathbb{C}_{(\mathbb{B}_0, \epsilon_1)}$};
  \node (b) at (0,4) {$\mathbb{C}_{(\mathbb{B}_0, \epsilon_2)}$};

  \node (d) at (3,0) {$\mathbb{C}_{(\mathbb{B}_1, \epsilon_0)}$};
  \node (e) at (6,0){$\mathbb{C}_{(\mathbb{B}_2, \epsilon_0)}$};

  \node (g) at (3,2) {$\mathbb{C}_{(\mathbb{B}_1, \epsilon_1)}$};
  \node (h) at (3,4) {$\mathbb{C}_{(\mathbb{B}_1, \epsilon_2)}$};

  \node (j) at (6,2) {$\mathbb{C}_{(\mathbb{B}_2, \epsilon_1)}$};
  \node (k) at (6,4) {$\mathbb{C}_{(\mathbb{B}_2, \epsilon_2)}$};

 \draw[->] (o) -- (a) ;
 \draw[->] (a) -- (b);
 
 \draw[->](o) -- (d);
 \draw[->](d) -- (e);

 \draw[->](a) -- (g);
 \draw[->](g) -- (j);

 \draw[->](d) -- (g);
 \draw[->](g) -- (h);

 \draw[->](b) -- (h);
 \draw[->](h) -- (k);

  \draw[->](e) -- (j);
 \draw[->](j) -- (k);

  \draw[dotted, thick] (6.9, 0) --(7.14, 0);
  \draw[dotted, thick] (6.9, 2) --(7.14, 2);
  \draw[dotted, thick] (6.9, 4) --(7.14, 4);

 \draw[dotted, thick] (6., 4.37) --(6, 4.6);
 \draw[dotted, thick] (3., 4.37) --(3, 4.6);
 \draw[dotted, thick] (0., 4.37) --(0, 4.6);

\end{tikzpicture}}
\hspace{0.15cm}
\resizebox{1.95in}{1.5in}{\begin{tikzpicture}
 \node (o) at (0,0) {$\mathbb{N}(\mathbb{C}_{(\mathbb{B}_0, \epsilon_0)})$};
  
  \node (a) at (0,2) {$\mathbb{N}(\mathbb{C}_{(\mathbb{B}_0, \epsilon_1)})$};
  \node (b) at (0,4) {$\mathbb{N}(\mathbb{C}_{(\mathbb{B}_0, \epsilon_2)})$};
  
  \node (d) at (3,0) {$\mathbb{N}(\mathbb{C}_{(\mathbb{B}_1, \epsilon_0)})$};
  \node (e) at (6,0){$\mathbb{N}(\mathbb{C}_{(\mathbb{B}_2, \epsilon_0)})$};

  \node (g) at (3,2) {$\mathbb{N}(\mathbb{C}_{(\mathbb{B}_1, \epsilon_1)})$};
  \node (h) at (3,4) {$\mathbb{N}(\mathbb{C}_{(\mathbb{B}_1, \epsilon_2)})$};

  \node (j) at (6,2) {$\mathbb{N}(\mathbb{C}_{(\mathbb{B}_2, \epsilon_1)})$};
  \node (k) at (6,4) {$\mathbb{N}(\mathbb{C}_{(\mathbb{B}_2, \epsilon_2)})$};

 \draw[->] (o) -- (a) ;
 \draw[->] (a) -- (b);

 \draw[->](o) -- (d);
 \draw[->](d) -- (e);

 \draw[->](a) -- (g);
 \draw[->](g) -- (j);

 \draw[->](d) -- (g);
 \draw[->](g) -- (h);

 \draw[->](b) -- (h);
 \draw[->](h) -- (k);

  \draw[->](e) -- (j);
 \draw[->](j) -- (k);

  \draw[dotted, thick] (7.2, 0) --(7.44, 0);
  \draw[dotted, thick] (7.2, 2) --(7.44, 2);
  \draw[dotted, thick] (7.2, 4) --(7.44, 4);

 \draw[dotted, thick] (6., 4.37) --(6, 4.6);
 \draw[dotted, thick] (3., 4.37) --(3, 4.6);
 \draw[dotted, thick] (0., 4.37) --(0, 4.6);
  
\end{tikzpicture}}
\hspace{0.1cm}
\resizebox{1.95in}{1.5in}{\begin{tikzpicture}
  \node (o) at (0,0) {$H_k(\mathbb{N}(\mathbb{C}_{(\mathbb{B}_0, \epsilon_0)}))$};
  
  \node (a) at (0,2) {$H_k(\mathbb{N}(\mathbb{C}_{(\mathbb{B}_0, \epsilon_1)}))$};
  \node (b) at (0,4) {$H_k(\mathbb{N}(\mathbb{C}_{(\mathbb{B}_0, \epsilon_2)}))$};
  
  \node (d) at (4,0) {$H_k(\mathbb{N}(\mathbb{C}_{(\mathbb{B}_1, \epsilon_0)}))$};
  \node (e) at (8,0){$H_k(\mathbb{N}(\mathbb{C}_{(\mathbb{B}_2, \epsilon_0)}))$};
  
  \node (g) at (4,2) {$H_k(\mathbb{N}(\mathbb{C}_{(\mathbb{B}_1, \epsilon_1)}))$};
  \node (h) at (4,4) {$H_k(\mathbb{N}(\mathbb{C}_{(\mathbb{B}_1, \epsilon_2)}))$};

  \node (j) at (8,2) {$H_k(\mathbb{N}(\mathbb{C}_{(\mathbb{B}_2, \epsilon_1)}))$};
  \node (k) at (8,4) {$H_k(\mathbb{N}(\mathbb{C}_{(\mathbb{B}_2, \epsilon_2)}))$};

 \draw[->] (o) -- (a) ;
 \draw[->] (a) -- (b);

 \draw[->](o) -- (d);
 \draw[->](d) -- (e);

 \draw[->](a) -- (g);
 \draw[->](g) -- (j);

 \draw[->](d) -- (g);
 \draw[->](g) -- (h);

 \draw[->](b) -- (h);
 \draw[->](h) -- (k);

  \draw[->](e) -- (j);
  \draw[->](j) -- (k);

  \draw[dotted, thick] (9.4, 0) --(9.67, 0);
  \draw[dotted, thick] (9.4, 2) --(9.67, 2);
  \draw[dotted, thick] (9.4, 4) --(9.67, 4);

 \draw[dotted, thick] (8, 4.37) --(8, 4.6);
 \draw[dotted, thick] (4., 4.37) --(4, 4.6);
 \draw[dotted, thick] (0., 4.37) --(0, 4.6);
  
\end{tikzpicture}}

\caption{Left to Right: Bi-filtration of cluster covers of $X$; Bi-filtration of the nerve of covers of $X$; Bi-filtration of the homology groups of the nerve of covers of $X$. $\mathbb{B}$ increases in the positive $x$-direction, and $\epsilon$ increases in the positive $y$-direction.}
\label{TDF_nc}
\end{figure}
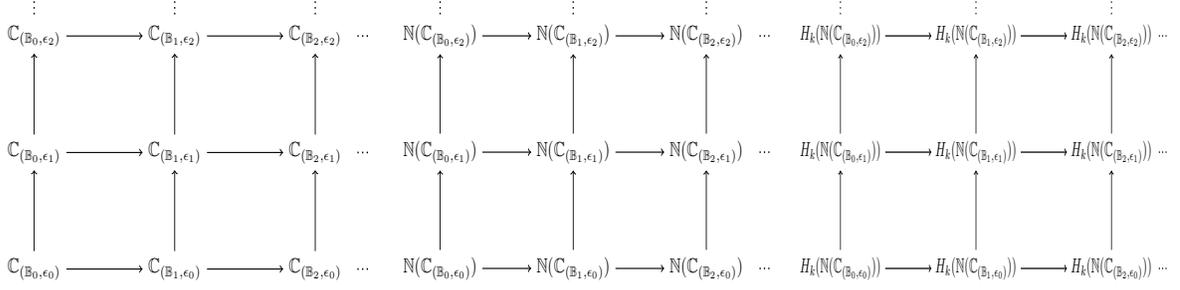

  We now turn our focus to a two dimensional filtration (bi-filtration) of covers of a dataset, $X$, and its $\delta$-perturbation, $X_{\delta}$. When DBSCAN is used to cluster dataset, $X$, we have shown in section \ref{dbsec} that there is a filtration of cluster covers as $\mathbb{B}$ increases, as $\epsilon$ increases, and as $MinPts$ decreases. A three dimensional filtration can be obtained. However in this section, we fix the parameter $MinPts$ and allow the parameters $\mathbb{B}$ and $\epsilon$ to vary. Hence, we only focus on a bi-filtration where one filtration is in the direction of $\mathbb{B}$ and the other filtration is in the direction of $\epsilon$. Note that a bi-filtration can be obtained with any of the two (out of the three) parameters. 
  Given a dataset, $X \in \mathbb{R}^n$, we have proven that there exists a filtration of cluster covers as $\mathbb{B}$ increases. So, we get a filtration of simplicial complexes (or mapper graphs) as $\mathbb{B}$ increases; that is a filtration of cluster covers 
$$\{c^{\mathbb{B}_i, \mathbb{B}_j} : \mathbb{C}_{\mathbb{B}_i} \rightarrow \mathbb{C}_{\mathbb{B}_j}, \forall \mathbb{B}_i \leq \mathbb{B}_j\}$$ 
induces a filtration of simplicial complexes 
$$\{\Phi^{\mathbb{B}_i, \mathbb{B}_j} : \mathbb{N}(\mathbb{C}_{\mathbb{B}_i}) \rightarrow\mathbb{N}(\mathbb{C}_{\mathbb{B}_j}), \forall \mathbb{B}_i \leq \mathbb{B}_j\}.$$
Let $k \in \mathbb{N}$ be fixed, then the above filtration of simplicial complexes induces a filtration of $k^{th}$ homology groups:
$$\{f^{\mathbb{B}_i, \mathbb{B}_j} : H_k(\mathbb{N}(\mathbb{C}_{\mathbb{B}_i})) \rightarrow H_k(\mathbb{N}(\mathbb{C}_{\mathbb{B}_j})), \forall \mathbb{B}_i \leq \mathbb{B}_j\}.$$
We have also proven that there exists a filtration of cluster covers as $\epsilon$ increases. Thus, we get a filtration of simplicial complexes as $\epsilon$ increases; that is a filtration of cluster covers 
$$\{c^{\epsilon_i, \epsilon_j} : \mathbb{C}_{\epsilon_i} \rightarrow \mathbb{C}_{\epsilon_j}, \forall \epsilon_i \leq \epsilon_j\}$$ 
induces a filtration of simplicial complexes 
$$\{\Phi^{\epsilon_i, \epsilon_j} : \mathbb{N}(\mathbb{C}_{\epsilon_i}) \rightarrow \mathbb{N}(\mathbb{C}_{\epsilon_j}), \forall \epsilon_i \leq \epsilon_j\}.$$
For a fixed $k \in \mathbb{N}$, the above filtration of simplicial complexes induces a filtration of $k^{th}$ homology groups:
$$\{f^{\epsilon_i, \epsilon_j} : H_k(\mathbb{N}(\mathbb{C}_{\epsilon_i})) \rightarrow H_k(\mathbb{N}(\mathbb{C}_{\epsilon_j})), \forall \epsilon_i \leq \epsilon_j\}.$$
Note that we have a bi-filtration of cluster covers, simplicial complexes, and homology groups, where $\mathbb{B}$ is the first dimension and $\epsilon$ is the second dimension. For any two $\mathbb{C}_{(\mathbb{B}_i, \epsilon_i)}$ and $\mathbb{C}_{(\mathbb{B}_j, \epsilon_j)}$, $(\mathbb{B}_i, \epsilon_i) \leq (\mathbb{B}_j, \epsilon_j)$ if and only if $\mathbb{B}_i \leq \mathbb{B}_j$ and $\epsilon_i \leq \epsilon_j$. We can visualize this bi-filtration of cluster covers, bi-filtration of simplicial complexes, and bi-filtration of homology groups as in figure \ref{TDF_nc}. Along any increasing path of these bi-filtrations of covers, simplicial complexes, and homology groups, we get filtrations presented as follows.

\begin{equation}\label{eqCdiag0}
\begin{split} 
\dots \rightarrow \mathbb{C}_{(\mathbb{B}_i, \epsilon_i)} \rightarrow \mathbb{C}_{(\mathbb{B}_j, \epsilon_j)} \rightarrow \mathbb{C}_{(\mathbb{B}_k, \epsilon_k)} \rightarrow \dots \\
\forall (\mathbb{B}_i, \epsilon_i) \leq (\mathbb{B}_j, \epsilon_j) \leq (\mathbb{B}_k, \epsilon_k). 
\end{split} 
\end{equation}

\begin{equation}\label{eqCdiag1}
\begin{split}
\dots \rightarrow \mathbb{N}(\mathbb{C}_{(\mathbb{B}_i, \epsilon_i)}) \rightarrow \mathbb{N}(\mathbb{C}_{(\mathbb{B}_j, \epsilon_j)}) \rightarrow \mathbb{N}(\mathbb{C}_{(\mathbb{B}_k, \epsilon_k)}) \rightarrow \dots \\
\forall (\mathbb{B}_i, \epsilon_i) \leq (\mathbb{B}_j, \epsilon_j) \leq (\mathbb{B}_k, \epsilon_k). 
\end{split}
\end{equation}

\begin{equation}\label{eqCdiag2}
\begin{split}
\dots \rightarrow H_k(\mathbb{N}(\mathbb{C}_{(\mathbb{B}_i, \epsilon_i)})) \rightarrow H_k(\mathbb{N}(\mathbb{C}_{(\mathbb{B}_j, \epsilon_j)})) \rightarrow H_k(\mathbb{N}(\mathbb{C}_{(\mathbb{B}_k, \epsilon_k)})) \rightarrow \dots \\
\forall (\mathbb{B}_i, \epsilon_i) \leq (\mathbb{B}_j, \epsilon_j) \leq (\mathbb{B}_k, \epsilon_k)
\end{split}
\end{equation}

\par

  Suppose also DBSCAN is applied to cluster $X_{\delta}$ wrt. $\epsilon$ and $MinPts$, and as $\mathbb{B}$ increases, similar bi-filtrations exist.
 
\noindent For any two $\mathbb{D}_{(\mathbb{B}_i, \epsilon_i)}$ and $\mathbb{D}_{(\mathbb{B}_j, \epsilon_j)}$, $(\mathbb{B}_i, \epsilon_i) \leq (\mathbb{B}_j, \epsilon_j)$ if and only if $\mathbb{B}_i \leq \mathbb{B}_j$ and $\epsilon_i \leq \epsilon_j$. Note that along these bi-filtrations of covers, simplicial complexes, and homology groups, we have the following 1-D filtrations:
 \begin{equation}\label{eqDdiag0}
\begin{split}
\dots \rightarrow \mathbb{D}_{(\mathbb{B}_i, \epsilon_i)} \rightarrow \mathbb{D}_{(\mathbb{B}_j, \epsilon_j)} \rightarrow \mathbb{D}_{(\mathbb{B}_k, \epsilon_k)} \rightarrow \dots \\
\forall (\mathbb{B}_i, \epsilon_i) \leq (\mathbb{B}_j, \epsilon_j) \leq (\mathbb{B}_k, \epsilon_k)
\end{split}
\end{equation}

 \begin{equation}\label{eqDdiag1}
 \begin{split}
\dots \rightarrow \mathbb{N}(\mathbb{D}_{(\mathbb{B}_i, \epsilon_i)}) \rightarrow \mathbb{N}(\mathbb{D}_{(\mathbb{B}_j, \epsilon_j)}) \rightarrow \mathbb{N}(\mathbb{D}_{(\mathbb{B}_k, \epsilon_k)}) \rightarrow \dots \\
\forall (\mathbb{B}_i, \epsilon_i) \leq (\mathbb{B}_j, \epsilon_j) \leq (\mathbb{B}_k, \epsilon_k)
\end{split}
\end{equation}

 \begin{equation}\label{eqDdiag2}
 \begin{split}
\dots \rightarrow H_k(\mathbb{N}(\mathbb{D}_{(\mathbb{B}_i, \epsilon_i)})) \rightarrow H_k(\mathbb{N}(\mathbb{D}_{(\mathbb{B}_j, \epsilon_j)})) \rightarrow H_k(\mathbb{N}(\mathbb{D}_{(\mathbb{B}_k, \epsilon_k)})) \rightarrow \dots \\
\forall (\mathbb{B}_i, \epsilon_i) \leq (\mathbb{B}_j, \epsilon_j) \leq (\mathbb{B}_k, \epsilon_k)
\end{split}
\end{equation}

  As the two parameters, $\mathbb{B}$ and $\epsilon$, increase, the topology of the underlying space of $X$ and $X_{\delta}$ evolves, and persistence homology is used to study this evolution of the topology of the underlying space for both $X$ and $X_{\delta}$. In section \ref{stabBiFilt}, we investigate conditions under which the bi-filtration of the cluster covers, simplicial complexes, and homology groups of the dataset, $X$, and the bi-filtration of the cluster covers, simplicial complexes, and homology groups of the perturbed dataset, $X_{\delta}$, are $\xi$-interleaved.

\subsection{Interleaving of Bi-filtrations}\label{interBiFit}
  In this section, we review the definition of interleaving of two bi-filtrations to prove our result of stability. Barcodes (or persistence diagrams) captures the evolution of the topology of the underlying space of a data $X$ in a one-dimensional filtration. Gunnar Carlsson et. al. proved that there are no barcodes (or persistence diagrams) for multi-dimensional persistence modules \cite{nobarformulti}. Bottleneck distance is defined on barcodes. Multi-dimension filtrations call for a generalization of a bottleneck distance. Chazal et. al. introduced $\epsilon$-interleaving distance between two persistence modules\cite{chazel}, and Michael Lesnick proposed a generalization of  bottleneck distance $\xi$-interleaving distance using a category theory framework \cite{DBLP}\cite{MichaelThesis}. Michael Lesnick defines an $\xi$-interleaving between two \textit{$n$-graded modules}. 

\begin{definition}{\textit{$n$-graded module \cite{DBLP}:}}\label{ML0}
Let $P_{n}$ be a polynomial ring in $n$ variables $x = \{x_1, x_2, \dots, x_n\}$. An $n$-graded module is a $P_{n}$-module M such that $M \simeq \oplus_{a \in \mathbb{R}^n} M_{a}$ and $x^{b}(M_a) \subset M_{a+b}$ for all $a \in \mathbb{R}^n, b \in [0,\infty)^n$ where $M_{a}$ is a vector space over some field $k$.
\noindent The action of $x^{b-a}$ gives rise to a linear map $\varphi: M_{a} \rightarrow M_{b}$ for all $a \leq b \in \mathbb{R}^{n}$.
\end{definition}

  Note that in the context of this paper, $n$-graded module is a two dimensional persistence module (bi-filtration of homology groups) where $x = \{\mathbb{B}, \epsilon\}$. The linear map $\varphi$ in definition \ref{ML0} is analogous to any composition of homomorphism (horizontal, vertical, or diagonal) in section \ref{bifilt}.
\par
\begin{definition}{\textit{Shift \cite{DBLP}:}}\label{ML1}
For M an $n$-graded module and $v \in \mathbb{R}^n$, $M(v)$ is the shifted module such that $M(v)_u = M_{v+u}$.
\end{definition}

\begin{definition}{\textit{Transition Morphism \cite{DBLP}:}}\label{ML2}
For M be an $n$-graded module, $\bar{\xi} =\{\xi, \xi, \dots, \xi\} \in \mathbb{R}_{+}^n$, and $M(\bar{\xi})$,
$$\varphi_M^{\bar{\xi}} : M \rightarrow M(\bar{\xi})$$
is the (diagonal) $\xi$-transition morphism such that $\varphi^{{\bar{\xi}}}_M(M_a) = \varphi_M (a + \bar{\xi})$.
\end{definition}
  We will use the notation $\xi$ instead of $\bar{\xi}$
\begin{definition}{\textit{$\xi$-Interleavings \cite{chazel} \cite{DBLP}:}}\label{ML3}
Let $\xi \geq 0$. Two $n$-modules M and N are $\xi$-interleaved if there are morphisms $f: M \rightarrow N(\xi)$ and $g:N \rightarrow M(\xi)$ such that 
$\varphi_{N}^{2\xi} = f(\xi) \circ g $ and $\varphi_{M}^{2\xi} = g(\xi) \circ f $.

\end{definition}
 
  Definition \ref{ML3} implies that if there are morphisms $f$ and $g$ between filtrations in equation \ref{eqCdiag0} and  equation \ref{eqDdiag0}, equation \ref{eqCdiag1} and equation \ref{eqDdiag1}, and equation \ref{eqCdiag2} and equation \ref{eqDdiag2} such that the conditions of definition \ref{ML3} are satisfied, then the filtrations are $\xi$-interleaved. 

\par
\subsection{Stability Against Perturbation}\label{stabBiFilt}
In this section, we study interleavings between two bi-filtrations where the first bi-filtration is with respect to a dataset $X$, and the second bi-filtration is with respect to $\delta$-perturbation of $X$, $X_{\delta}$. We will first investigate the conditions under which a well-defined family of maps, $\phi$ and $\psi$, between a bi-filtration of the cluster covers of $X$ and a bi-filtration of the cluster covers of $X_{\delta}$ are obtained. Recall that covers of $X$ (and $X_{\delta}$) refer to clusters of $X$ (and $X_{\delta}$). When DBSCAN is used to cluster both $X$ and $X_{\delta}$, a relationship between the parameters ($\mathbb{B}$ and $\epsilon$) of $X$ and $X_{\delta}$ needs to be established in order to obtain a well-defined family of morphisms as in definition \ref{ML3}.
\par 
  We state the results of this section in proposition \ref{prop22} and corollaries \ref{cor2deltainter} and \ref{cor4deltainter}. Proposition \ref{prop22} proves the existence of morphisms between bi-filtration of cluster covers of $X$ and $X_{\delta}$ that satisfy the conditions of definition \ref{ML3}, and corollary \ref{cor2deltainter} shows that the two covers are $2\delta$-interleaved. Note $2\delta$-interleaved covers induce $2\delta$-interleaved simplicial complexes which in turn induces $2\delta$-interleaved bi-graded modules  \cite{deybook}\cite{dey}. Corollary \ref{cor4deltainter} proves similar result for bi-filtrations of homology groups of $X$ and $X_{\delta}$; that is, the two bi-filtrations of homology groups are $2\delta$-interleaved.

  Let us refer to the example in figure \ref{stblty12} and state the issues that arise when defining a well-defined map. First, there are points in $X_{\delta}$ not in $bin$ (figure \ref{stblty12}.B), hence no map exists from clusters of $X$ in $bin$ to clusters of $X_{\delta}$ in $bin$. A straight forward solution is to increase the size of $bin$ by $\delta$; denoted by $bin_ \delta\in \mathbb{B}_{\delta}$.

  Furthermore, if the distance between two adjacent points, $x$ and $y$ in $X$ (i.e. in $\mathbb{B}$) is $dist(x,y) = d$, then the distance between the perturbed points $x_{\delta}$ and $y_{\delta}$ in $X_{\delta}$ (i.e. in $\mathbb{B}_{\delta}$) is $dist(x,y) \leq d + 2\delta$. As shown in figure \ref{stblty12}.C, there are at four clusters in $bin_{\delta}$ wrt. $\epsilon=d$ and $MinPts=3$. Hence the map $\phi$ is not well-defined since $\phi$ sends one cluster of $X$ in $bin\in\mathbb{B}$ to three clusters of $X_{\delta}$ in $bin_{\delta}\in\mathbb{B}_{\delta}$. If DBSCAN is used to cluster $X_{\delta}$ wrt. $\epsilon= d + 2\delta$ and $MinPts = 3$, then there is a cluster, $C_1^{bin_{\delta}}$, in $bin_{\delta}$ wrt. $\epsilon = d + 2\delta$ and $MinPts=3$ that contains $C_1^{bin}$. There is also a cluster, $C_2^{bin_{\delta}}$, that contains $C_2^{bin}$. The containment of the clusters described here ought to be understood as follows: $x \in C_i^{bin} \implies x_{\delta} \in C_i^{bin_{\delta}}$. We now define $\phi : \mathbb{C}_{\mathbb{B}, \epsilon} \rightarrow \mathbb{D}_{\mathbb{B}_{\delta}, \epsilon+2\delta}$ where $\mathbb{C}_{\mathbb{B}, \epsilon} $ is a cover of $X$ obtained by applying DBSCAN in $\mathbb{B}$ wrt. $\epsilon$ and $MinPts$, and  $\mathbb{D}_{\mathbb{B}_{\delta}, \epsilon+2\delta}$ is a cover of $X_{\delta}$ obtained by applying DBSCAN in $\mathbb{B}_{\delta}$ wrt. $\epsilon+2\delta$ and $MinPts$. Hence, $\phi$ is well-defined.

    A similar argument is used to show that $\psi : \mathbb{D}_{\mathbb{B}_{\delta}, \epsilon+2\delta} \rightarrow \mathbb{C}_{\mathbb{B}_{2\delta}, \epsilon+4\delta}$ is well-defined where $\mathbb{C}_{\mathbb{B}_{2\delta}, \epsilon+4\delta}$ is a cover of $X$ obtained by applying DBSCAN in $\mathbb{B}_{2\delta}$ wrt. $\epsilon+4\delta$ and $MinPts$, and  $\mathbb{D}_{\mathbb{B}_{\delta}, \epsilon+2\delta}$ is a cover of $X_{\delta}$ obtained by applying DBSCAN in $\mathbb{B}_{\delta}$ wrt. $\epsilon+2\delta$ and $MinPts$. Proposition \ref{prop22} generalizes the existence of well-defined families of maps $\phi$ and $\psi$.

\begin{prop}\label{prop22}
Let $X$ and $X_{\delta}$ be datasets such that $X_{\delta}$ is obtained by perturbing $X$ by at most $\delta$. Let $k, l \in \mathbb{R}$. 
\vspace{0.1cm}
  Let $\mathcal{C} = \{c^{(k,l),(k+1, l+2)} : \mathbb{C}_{\mathbb{B}_{k\delta}, \epsilon+l\delta} \rightarrow \mathbb{C}_{\mathbb{B}_{(k+1)\delta}, \epsilon+(l+2)\delta},$ $\forall (\mathbb{B}_{k\delta}, \epsilon+l\delta)$ be a filtration of cluster covers of $X$ where $\mathbb{C}_{\mathbb{B}_{a}, \epsilon+b}$ is a cover of $X$ obtained by applying DBSCAN in $\mathbb{B}_{a}$ wrt. $\epsilon+b$ and $MinPts$
(hence, $\mathbb{C}_{\mathbb{B}_{k\delta}, \epsilon+l\delta}$ is a cover of $X$ obtained by applying DBSCAN in $\mathbb{B}_{k\delta}$ wrt. $\epsilon+l\delta$ and $MinPts$). 
\vspace{0.1cm}
  Let $\mathcal{D} = \{d^{(k,l),(k+1, l+2)} : \mathbb{D}_{\mathbb{B}_{k\delta}, \epsilon+l\delta} \rightarrow \mathbb{D}_{\mathbb{B}_{k\delta}, \epsilon+l\delta},$ $\forall (\mathbb{B}_{k\delta}, \epsilon+l\delta)$ be a filtration of cluster covers of $X_{\delta}$ where $\mathbb{D}_{\mathbb{B}_{a}, \epsilon+b}$ is a cover of $X_{\delta}$ obtained by applying DBSCAN in $\mathbb{B}_{a}$ wrt. $\epsilon+b$ and $MinPts$
(hence, $\mathbb{D}_{\mathbb{B}_{k\delta}, \epsilon+l\delta}$ is a cover of $X_{\delta}$ obtained by applying DBSCAN in $\mathbb{B}_{k\delta}$ wrt. $\epsilon+l\delta$ and $MinPts$). 
\vspace{0.1cm}
  Assume no free-border points. Then there are families of maps $\phi$ and $\psi$ such that the diagram in figure \ref{coverInterleave} commutes. 
\end{prop}
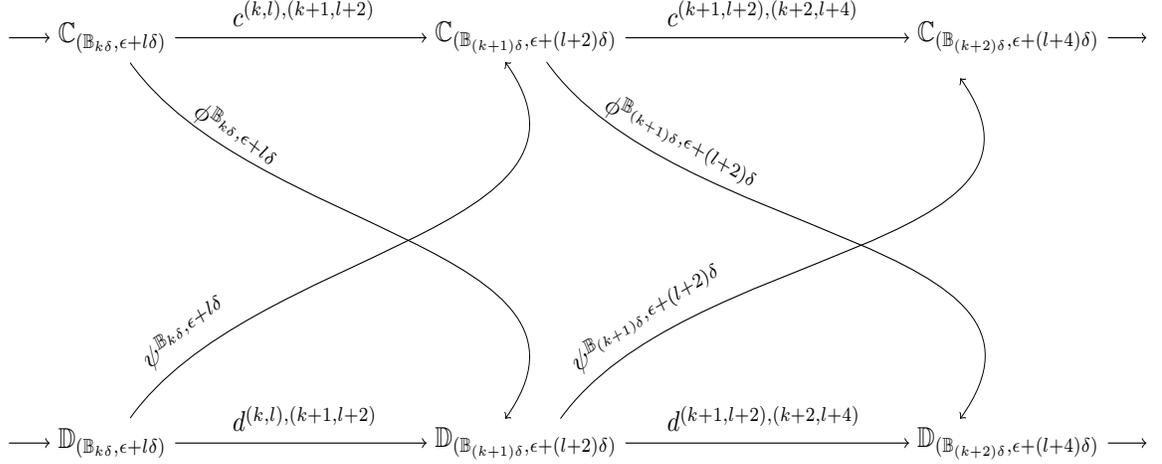
\begin{figure}[h!]
\centering
\resizebox{6in}{2.5in}{
\begin{tikzpicture}

\node (d) at (0,-1) {$\mathbb{D}_{(\mathbb{B}_{k\delta}, \epsilon + l\delta)}$};
\node (e) at (6,-1){$\mathbb{D}_{(\mathbb{B}_{(k+1)\delta}, \epsilon+ (l+2)\delta)}$};
\node (f) at (13,-1) {$\mathbb{D}_{(\mathbb{B}_{(k+2)\delta}, \epsilon + (l+4)\delta)}$};

\node (a) at (0,4) {$\mathbb{C}_{(\mathbb{B}_{k\delta}, \epsilon + l\delta)}$};
\node (b) at (6,4){$\mathbb{C}_{(\mathbb{B}_{(k+1)\delta}, \epsilon+ (l+2)\delta)}$};
\node (c) at (13,4) {$\mathbb{C}_{(\mathbb{B}_{(k+2)\delta}, \epsilon + (l+4)\delta)}$};

\node[label={[label distance=0.5cm,rotate=-35]right: $\phi^{\mathbb{B}_{k\delta}, \epsilon+l\delta}$ }] at (0.6, 3.5) {};
\node[label={[label distance=0.5cm,rotate=30]right: $\psi^{\mathbb{B}_{k\delta}, \epsilon + l\delta}$ }] at (-0.2, -0.3) {};
\node[label={[label distance=0.5cm,rotate=-27]right: $\phi^{\mathbb{B}_{(k+1)\delta}, \epsilon +(l+ 2)\delta}$ }] at (6.5, 3.5) {};
\node[label={[label distance=0.5cm,rotate=27]right: $\psi^{\mathbb{B}_{(k+1)\delta}, \epsilon+(l+2)\delta}$ }] at (6., -0.3) {};

\node[label={[label distance=0.5cm]right: $c^{(k,l),(k+1, l+2)}$ }] at (1., 4.3) {};
\node[label={[label distance=0.5cm]right: $c^{(k+1,l+2),(k+2, l+4)}$ }] at (7.3, 4.3) {};
\node[label={[label distance=0.5cm]right: $d^{(k,l),(k+1, l+2)}$ }] at (1., -0.7) {};
\node[label={[label distance=0.5cm]right: $d^{(k+1,l+2),(k+2, l+4)}$ }] at (7.3, -0.7) {};

\draw[->](-1.5,4) -- (a);  
\draw[->](a) -- (b);
\draw[->](b) -- (c);
\draw[->](c) -- (15,4);
\draw[->](-1.5,-1) -- (d);  
\draw[->](d) -- (e);
\draw[->](e) -- (f);
\draw[->](f) -- (15,-1);
  
\draw[->](a) to  [out=-50, in=50] (5.7,-0.7);
\draw[->](6.5,-0.7) to [out=50, in=-50] (12.3,3.5);
\draw[->](6.3,3.7) to [out=-50, in=50] (12.3, -0.7);
\draw[->](d) to  [out=50, in=-50](5.7,3.7);
     
     \end{tikzpicture}}
\caption{Morphisms between the bi-filtrations of the covers of $X$ and $X_{\delta}$.}
\label{coverInterleave}  
\end{figure}
\begin{proof}

Let us first define $\phi$ as follows. Let $C \in \mathbb{C}_{\mathbb{B}_{k\delta}, \epsilon+l\delta}$ be a cluster wrt. $\epsilon+l\delta$ and $MinPts$ in $\mathbb{B}_{k\delta}$, and $\phi^{\mathbb{B}_{k\delta}, \epsilon+l\delta}$ is defined on the clusters in $\mathbb{C}_{\mathbb{B}_{k\delta}, \epsilon+l\delta}$. Suppose two core points $p$ and $q$ represent $C$. Then $p$ is density-reachable from $q$ wrt. $\epsilon +l\delta$ and $MinPts$ in $\mathbb{B}_{k\delta}$ (and vice-versa). That is, there is a sequence of core points  $\{p = q_1, q_2, \dots, q_n=q\}$ such that $dist(q_i, q_{i+1}) \leq  \epsilon+l\delta$. That is, $q_i$ is directly density-reachable from $q_{i+1}$ and visa versa wrt. $\epsilon +l\delta$ and $MinPts$. Moreover, if $x \in C$, then $dist(x, q_i)\leq \epsilon+l\delta$ for some $i\in\{1,2,3, \dots, n\}$.
 
  \textit{Claim:} There is a cluster $D \in \mathbb{D}_{\mathbb{B}_{(k+1)\delta}, \epsilon+(l+2)\delta}$ wrt. $\epsilon+(l+2)\delta$ and $MinPts$ in $\mathbb{B}_{(k+1)\delta}$ represented by $p_{\delta}$ and $q_{\delta}$ where $dist(p, p_{\delta})\leq\delta$ and $dist(q, q_{\delta})\leq\delta$ such that if $x\in C$, $x_{\delta}\in D$.

  \textit{Proof of Claim:} Note $dist(q_{i\delta}, q_{(i+1)\delta}) \leq  \epsilon+(l+2)\delta$, then $q_{i\delta}$ is density-reachable from $q_{(i+1)\delta}$ and vice versa wrt. $\epsilon+(l+2)\delta$ and $MinPts$ in $\mathbb{B}_{(k+1)\delta}$. For a point $x\in C$, $dist(x_{\delta}, q_{i\delta})\leq \epsilon+(l+2)\delta$ for some $i\in\{1,2,3, \dots, n\}$ in $\mathbb{B}_{(k+1)\delta}$. It follows $p_{\delta}$ and $q_{\delta}$ determine the cluster $D$ in $\mathbb{D}_{\mathbb{B}(k+1)\delta, \epsilon +(l+2)\delta}$ wrt. $\epsilon+(l+2)\delta$ and $MinPts$. Thus the claim is proved.
\vspace{0.25cm}
  Define $\phi^{\mathbb{B}_{k\delta}, \epsilon +l\delta}(C) = D$; hence a cluster $C \in \mathbb{C}_{\mathbb{B}_{k\delta}, \epsilon+l\delta}$ is mapped to a cluster $D = \phi^{\mathbb{B}_{k\delta}, \epsilon +l\delta}(C) \in \mathbb{D}_{\mathbb{B}_{(k+1)\delta}, \epsilon+(l+2)\delta}$ such that if $x$ is a point in $C$ then $x_{\delta}$ is in $D$. Since the definition of $\psi^{\mathbb{B}_{k\delta}, \epsilon +l\delta}$ follows similar argument to that of $\phi^{\mathbb{B}_{k\delta}, \epsilon +l\delta}$, $\psi^{\mathbb{B}_{k\delta}, \epsilon +l\delta}$ is also well-defined. Moreover, both $\phi^{\mathbb{B}_{k\delta}, \epsilon +l\delta}$ and $\psi^{\mathbb{B}_{k\delta}, \epsilon +l\delta}$ are inclusion maps in the sense that they both map a cluster to its perturbation.
   Now let us show the commutativity of the diagram in figure \ref{coverInterleave}. It is fairly straight forward to see the commutativity of the diagram since all the maps are inclusion maps. Therefore,
\begin{itemize}

\item $\phi^{\mathbb{B}_{(k+1)\delta}, \epsilon + (l+2)\delta} \circ c^{(k,l),(k+1, l+2)} = d^{(k+1,l+2),(k+3,l+4)} \circ \phi^{\mathbb{B}_{k\delta}, \epsilon + l\delta}$
\item $\psi^{\mathbb{B}_{(k+1)\delta}, \epsilon + (l+2)\delta} \circ d^{(k,l),(k+1, l+2)} = c^{(k+1,l+2),(k+3,l+4)} \circ \psi^{\mathbb{B}_{k\delta}, \epsilon + l\delta}$

\item $c^{(k,l),(k+3, l+6)} = \psi^{\mathbb{B}_{(k+2)\delta}, \epsilon + (l+4)\delta} \circ d^{(k+1,l+2),(k+3,l+4)} \circ \phi^{\mathbb{B}_{k\delta}, \epsilon + l\delta} $
\item $d^{(k,l),(k+3, l+6)} = \phi^{\mathbb{B}_{(k+2)\delta}, \epsilon + (l+4)\delta} \circ c^{(k+1,l+2),(k+3,l+4)} \circ \psi^{\mathbb{B}_{k\delta}, \epsilon + l\delta} $

\end{itemize}

\end{proof}

\begin{corollary}\label{cor2deltainter}
Let $X$ and $X_{\delta}$ be datasets such that $X_{\delta}$ is obtained by perturbing $X$ by $\delta$. Let DBSCAN be applied to both $X$ and $X_{\delta}$, and assume no free-border point exists. Then, the bi-filtration of cluster covers of $X$ and the bi-filtration of cluster covers of $X_{\delta}$ are $2\delta$-interleaved. 

\end{corollary}

\begin{proof}
Let $X$ and $X_{\delta}$ be datasets such that $X_{\delta}$ is obtained by perturbing $X$ by $\delta$. Let $k, l \in \mathbb{R}$. 
\par
\noindent Let $\mathcal{C} = \{c^{(k,l),(k+2, l+2)} : \mathbb{C}_{\mathbb{B}_{k\delta}, \epsilon+l\delta} \rightarrow \mathbb{C}_{\mathbb{B}_{(k+2)\delta}, \epsilon+(l+2)\delta},$ $\forall (\mathbb{B}_{k\delta}, \epsilon+l\delta)$ be a filtration of cluster covers of $X$ such that $\mathbb{C}_{\mathbb{B}_{k\delta}, \epsilon+l\delta}$ is a cover of $X$ obtained by applying DBSCAN in $\mathbb{B}_{k\delta}$ wrt. $\epsilon+l\delta$ and $MinPts$. 
\par
\noindent Let $\mathcal{D} = \{d^{(k,l),(k+2, l+2)} : \mathbb{D}_{\mathbb{B}_{k\delta}, \epsilon+l\delta} \rightarrow \mathbb{D}_{\mathbb{B}_{k\delta}, \epsilon+l\delta}, \forall (\mathbb{B}_{k\delta}, \epsilon+l\delta) \leq (\mathbb{B}_{(k+2)\delta}, \epsilon+(l+2)\delta)\}$ be a filtration of cluster covers of $X_{\delta}$ such that $\mathbb{D}_{\mathbb{B}_{k\delta}, \epsilon+l\delta}$ is a cover of $X_{\delta}$ obtained by applying DBSCAN in $\mathbb{B}_{k\delta}$ wrt. $\epsilon+l\delta$ and $MinPts$. 
\par
\noindent Then, by proposition \ref{prop22}, there are families of maps $\phi$ and $\psi$ such that the diagram in figure \ref{coverInterleave1} commutes. And, conditions of definition \ref{ML3} are satisfied as shown in the proof of proposition \ref{prop22}. By identify the morphisms in definition \ref{ML3} to the morphisms in figure \ref{coverInterleave1}, we obtain that the bi-filtration of cluster covers of $X$ and the bi-filtration of cluster covers of $X_{\delta}$ are $2\delta$-interleaved. 
$$c^{(k,l), (k+4, l+4)}=\varphi_{N}^{4\delta}  \thinspace (Transition \thinspace morphism)$$
$$\psi^{\mathbb{B}_{(k+2)\delta}, \epsilon + (l+2)\delta} = f(2\delta)$$
$$\phi^{\mathbb{B}_{k\delta}, \epsilon+l\delta} = g$$
$$d^{(k,l), (k+4, l+4)} = \varphi_{M}^{4\delta}$$
$$\phi^{\mathbb{B}_{(k+2)\delta}, \epsilon + (l+2)\delta} = g(2\delta)$$
$$ \psi^{\mathbb{B}_{k\delta}, \epsilon+l\delta} = f $$

\end{proof}

\begin{figure}
\centering
\resizebox{5.3in}{2in}{
\begin{tikzpicture}

\node (d) at (0,-1) {$H_k(\mathbb{N}(\mathbb{D}_{(\mathbb{B}_k\delta, \epsilon + l\delta)}))$};
\node (e) at (6,-1){$H_k(\mathbb{N}(\mathbb{D}_{(\mathbb{B}_{(k+2)\delta}, \epsilon+ (l+2)\delta)}))$};
\node (f) at (13,-1) {$H_k(\mathbb{N}(\mathbb{D}_{(\mathbb{B}_{(k+2)\delta}, \epsilon + (l+4)\delta)}))$};

\node (a) at (0,4) {$H_k(\mathbb{N}(\mathbb{C}_{(\mathbb{B}_k\delta, \epsilon + l\delta)}))$};
\node (b) at (6,4){$H_k(\mathbb{N}(\mathbb{C}_{(\mathbb{B}_{(k+2)\delta}, \epsilon+ (l+2)\delta)}))$};
\node (c) at (13,4) {$H_k(\mathbb{N}(\mathbb{C}_{(\mathbb{B}_{(k+2)\delta}, \epsilon + (l+4)\delta)}))$};

\node[label={[label distance=0.5cm,rotate=-35]right: $\phi^{\mathbb{B}_{k\delta}, \epsilon+l\delta}$ }] at (0.6, 3.5) {};
\node[label={[label distance=0.5cm,rotate=30]right: $\psi^{\mathbb{B}_{k\delta}, \epsilon + l\delta}$ }] at (-0.2, -0.3) {};
\node[label={[label distance=0.5cm,rotate=-27]right: $\phi^{\mathbb{B}_{(k+2)\delta}, \epsilon +(l+ 2)\delta}$ }] at (6.5, 3.5) {};
\node[label={[label distance=0.5cm,rotate=27]right: $\psi^{\mathbb{B}_{(k+2)\delta}, \epsilon+(l+2)\delta}$ }] at (6., -0.3) {};

\node[label={[label distance=0.5cm]right: $c^{(k,l),(k+2, l+2)}$ }] at (1., 4.3) {};
\node[label={[label distance=0.5cm]right: $c^{(k+2,l+2),(k+2, l+4)}$ }] at (7.3, 4.3) {};
\node[label={[label distance=0.5cm]right: $d^{(k,l),(k+2, l+2)}$ }] at (1., -0.7) {};
\node[label={[label distance=0.5cm]right: $d^{(k+2,l+2),(k+2, l+4)}$ }] at (7.3, -0.7) {};

\draw[->](-1.5,4) -- (a);  
\draw[->](a) -- (b);
\draw[->](b) -- (c);
\draw[->](c) -- (15,4);
\draw[->](-1.5,-1) -- (d);  
\draw[->](d) -- (e);
\draw[->](e) -- (f);
\draw[->](f) -- (15,-1);
  
\draw[->](a) to  [out=-50, in=50] (5.7,-0.7);
\draw[->](6.5,-0.7) to [out=50, in=-50] (12.3,3.5);
\draw[->](6.3,3.7) to [out=-50, in=50] (12.3, -0.7);
\draw[->](d) to  [out=50, in=-50](5.7,3.7);

   \end{tikzpicture}}
\caption{The bi-filtrations of the covers of $X$ and $X_{\delta}$ are $2\delta$-interleaved.}
\label{coverInterleave1}  
\end{figure}

  The families of maps in proposition \ref{prop22} are inclusion maps, hence the bi-filtrations of the simplicial complexes of $X$ and $X_{\delta}$ are $2\delta$-interleaved. That is, the families of maps in proposition \ref{prop22} induce a a family of contiguous simplicial maps between the bi-filtrations. Two contiguous simplicial maps induce two equal homomorphism on the homology groups, hence the bi-filtrations of the homology groups of $X$ and $X_{\delta}$ are $2\delta$-interleaved. That is, the families of contiguous simplicial maps induce a family of homomorphisms between the bi-filtrations the homology groups  \cite{deybook}\cite{dey}. 
\begin{corollary}\label{cor4deltainter}
Let $X$ and $X_{\delta}$ be datasets such that $X_{\delta}$ is obtained by perturbing $X$ by $\delta$. Let DBSCAN be used to cluster both $X$ and $X_{\delta}$, and assume no free-border point exists. Then, the bi-filtration of the homology of $X$ and the bi-filtration of the homology of $X_{\delta}$ are $2\delta$-interleaved. 

\end{corollary}


\end{document}